\documentclass[11pt]{article}
\usepackage[T1]{fontenc}
\usepackage{lmodern}
\usepackage[margin=1in]{geometry}
\usepackage{amsmath,amssymb,amsthm}
\usepackage{enumitem,microtype}
\usepackage{xcolor}
\usepackage[colorlinks=true,linkcolor=blue!50!black,citecolor=blue!50!black,urlcolor=blue!50!black]{hyperref}
\hypersetup{pdftitle={Model Theory and Algorithms for One-Relator Groups},pdfauthor={Olga Kharlampovich and Alexei Myasnikov}}
\newtheorem{theorem}{Theorem}[section]
\newtheorem{lemma}[theorem]{Lemma}
\newtheorem{proposition}[theorem]{Proposition}
\newtheorem{corollary}[theorem]{Corollary}
\theoremstyle{definition}
\newtheorem{definition}[theorem]{Definition}
\newtheorem{citedresult}[theorem]{Cited result}

\newtheorem{example}[theorem]{Example}
\theoremstyle{remark}
\newtheorem{remark}[theorem]{Remark}
\newcommand{\BS}{\operatorname{BS}}
\newcommand{\Z}{\mathbb Z}
\newcommand{\Q}{\mathbb Q}

\newcommand{\Th}{\operatorname{Th}}
\newcommand{\CT}{\mathsf{CT}}
\newcommand{\Mal}{\mathsf{Mal}}
\newcommand{\Sq}{\mathsf{Sq}}
\newcommand{\CSA}{\mathsf{CSA}}
\newcommand{\Hyp}{\mathsf{Hyp}}
\newcommand{\length}{\ell}

\newcommand{\Comm}{\operatorname{Comm}}
\newcommand{\cdim}{\operatorname{cd}}

\newcommand{\Div}{\operatorname{Div}}

\setlist[enumerate]{leftmargin=*,itemsep=4pt,topsep=5pt}
\setlist[itemize]{leftmargin=*,itemsep=3pt,topsep=5pt}
\title{Model Theory and Algorithms\\for One-Relator Groups}
\author{Olga Kharlampovich \and Alexei Miasnikov}
\date{}

\newcommand{\Aut}{\operatorname{Aut}}
\newcommand{\cdq}{\operatorname{cd}_{\Q}}
\newcommand{\Fix}{\operatorname{Fix}}
\newcommand{\Inn}{\operatorname{Inn}}
\newcommand{\rank}{\operatorname{rank}}
\newcommand{\PP}{\mathcal P}

\begin{document}
\maketitle

\begin{abstract}
We study the elementary theories and algorithmic properties of one-relator
groups. We give a first-order sentence characterizing Baumslag--Solitar-freeness
in this class. For fixed $|m|>1$, we regularly interpret $\BS(1,m)$ as a
definable subgroup in every one-relator group containing it. We uniformly,
regularly and injectively interpret integer arithmetic in countable groups
of cohomological dimension at most two.  We obtain one uniform interpretation for all
one-relator groups containing a non-unimodular Baumslag--Solitar subgroup.
Every fixed finitely presented group is regularly interpretable in these
groups, and their elementary theories are undecidable. We also describe
the GBS overgroups of a non-unimodular Baumslag--Solitar subgroup.
For the unimodular groups we construct finite-index direct products of
free groups and $\Z$, and reduce the elementary theory of $\BS(k,k)$ to
that of a free group with a named periodic automorphism. 

  Using the recent Open AI result that
BS-free one-relator groups are hyperbolic we obtain a terminating algorithm recognizing this class and derive consequences for isomorphism, conjugacy, equations, and elementary theories. From the virtual compact specialness result of another Open AI paper, we obtain an effective construction of virtual free-by-cyclic presentations for hyperbolic one-relator groups, the free kernel may have infinite rank.
We prove that a finitely generated torsion-free virtually free-by-cyclic
group is toral relatively hyperbolic if and only if it is CSA, equivalently
if it contains neither $F_2\times\Z$ nor the Klein bottle group. For classes
of torsion-free one-relator groups in which CSA implies toral relative
hyperbolicity, we give a uniform algorithm deciding CSA and, in the
affirmative case, the existence of a subgroup isomorphic to $\Z^2$. The  model-theoretic results above are independent of the OpenAI results. \end{abstract}

\section{Introduction}

A one-relator group has a finite presentation
\[
G=\langle X\mid r\rangle.
\]
One-relator groups with torsion are hyperbolic \cite{Newman}. They are
residually finite by Wise's theorem \cite{Wise}, and hereditarily conjugacy
separable by \cite{MZ}. Their isomorphism problem is decidable by the
isomorphism algorithm for hyperbolic groups \cite{DG}.
Torsion-free one-relator groups have more varied properties. For example,
\(\BS(2,3)\) is not residually finite.

\begin{definition}\label{def:bs}
For nonzero integers \(m,n\), the Baumslag--Solitar group is
\[
\BS(m,n)=\langle a,t\mid ta^mt^{-1}=a^n\rangle.
\]
A group is \emph{BS-free} if it contains no subgroup isomorphic to any
\(\BS(m,n)\). In particular, \(\BS(1,1)=\Z^2\) is excluded.
We call \(\BS(m,n)\) \emph{unimodular} if \(|m|=|n|\), and
\emph{non-unimodular} otherwise.
\end{definition}

Throughout the paper, \([x,y]=xyx^{-1}y^{-1}\), and \(A^d=d^{-1}Ad\).
The elementary theory \(\Th(G)\) is the set of first-order sentences in
the language of groups which hold in \(G\). The Diophantine problem asks
whether a finite system of equations with coefficients in \(G\) has a
solution in \(G\).

Our first result, Theorem~\ref{thm:fo}, gives a finite first-order
characterization of BS-freeness among torsion-free one-relator groups.
Corollary~\ref{cor:all-fo} extends it to groups with torsion. Thus a
BS-free one-relator group cannot be elementarily equivalent to a
one-relator group containing a Baumslag--Solitar subgroup.

In Section~\ref{sec:interpretation}, we prove that, for fixed \(|m|>1\),
the group \(\BS(1,m)\) is regularly interpretable in every one-relator
group $G$ containing  it (Theorem~\ref{thm:bs-interpretation}). The admissible
parameters are the pairs satisfying \(a\neq1\) and \(tat^{-1}=a^m\).
The interpretation gives a definable subgroup isomorphic to \(\BS(1,m)\),
with the inherited group operations. Together with the interpretation of
arithmetic $(\Z,+,\mid,0,1)$ in solvable Baumslag--Solitar groups \cite{Khelif,DMBS}, this
implies undecidability of the elementary theory of $G$.

Section~\ref{sec:nonuni-bs} proves a general arithmetic theorem.
Let \(U\) be a countable group of integral cohomological dimension at most two, and let \(a\in U\) have infinite order. Suppose that \(\langle a\rangle\) is commensurated by \(U\).
If
\[
 sa^ms^{-1}=a^n,\qquad |m|\ne|n|,
\]
then arithmetic is uniformly, regularly and injectively interpretable
in $U$ (Theorem~\ref{arith:thm:criterion}). The formulas are independent
of $m,n$. The proof uses centers of nonabelian centralizers, cyclic
centralizers of elements $s$ with $|\Delta(s)|\ne1$, and
intersections of cyclic subgroups with their conjugates. Divisibility
of exponents is defined by inclusions between subgroups consisting of
commutators. The classical definability of multiplication from addition
and divisibility completes the interpretation. This theorem applies
to non-unimodular GBS groups and to their increasing unions of
cohomological dimension two.

Section~\ref{sec:nonuni-ambient} applies this result to one-relator
groups. If $H=\langle a,s\rangle\cong\BS(m,n)$ is non-unimodular, its
commensurator $M=\Comm_G(\langle a\rangle)$ is definable and is an
increasing union of non-unimodular GBS overgroups of $H$.
Every GBS overgroup of $H$ lies in $M$. In particular, an
inclusion-maximal such overgroup exists if and only if $M$ is finitely
generated, and then it is $M$.
We give a definable nonempty family of parameters for the intermediate
subgroups and obtain a uniform regular injective interpretation of
arithmetic in all one-relator groups containing a non-unimodular
Baumslag--Solitar subgroup (Theorem~\ref{ambient:thm:main}).

Every fixed finitely presented group is therefore uniformly, regularly
and injectively interpretable in these groups by arithmetic coding.
The resulting group operation need not be the inherited group operation.
One computable translation reduces integer arithmetic to the elementary
theory of every group in the class. These results concern the full
elementary theory; the Diophantine problems are considered in
Section~\ref{sec:beyond}.

We also discuss bi-interpretability. The soluble groups in
Theorem~\ref{arith:thm:criterion} are precisely $\BS(1,k)$, $|k|>1$.
The known result for $k>1$ \cite{DMBS}, together with a definable
reconstruction from an index-two subgroup for $k<-1$, gives
bi-interpretability with $\Z$. An example of cohomological dimension
two shows that this conclusion does not hold for the whole countable
class.

The unimodular groups \(\BS(k,k)\) and \(\BS(k,-k)\), \(k>1\), contain
explicit finite-index direct products of a free group and \(\Z\).
Their direct-product subgroups have decidable elementary theory.
Section~\ref{sec:unimodular} also describes the cyclic action giving
\(\BS(k,k)=F_k\rtimes\Z\) and proves a conditional decidability result.
The remaining requirement is decidability of the free-group theory with
the basis permutation named as a function.

For the recognition and algorithmic results we use the following result,
stated as Theorem~1.1 in \cite{OAI} and denoted below by (H).

\begin{citedresult}[H]\label{ass:H}
Every finitely generated BS-free one-relator group is word-hyperbolic.
\end{citedresult}

We also use the virtual compact specialness result stated as Theorem~1.1
in \cite{OAISpecial}. We have read both papers, but have not verified all
details of their proofs. Our conclusions that use these results are
conditional on their validity, and we indicate these dependencies
explicitly.

All our original model-theoretic results are independent of both cited
results. These include the first-order characterization of BS-freeness
in Section~\ref{sec:fo}, the interpretation and undecidability results in
Sections~\ref{sec:interpretation}--\ref{sec:nonuni-ambient}, and the
model-theoretic results on unimodular groups in Section~\ref{sec:unimodular}.
They originate in our earlier manuscript, written in 2019.
The applications of known results on hyperbolic groups to all BS-free
one-relator groups, including Corollary~\ref{cor:model}, use (H).

Under (H), Theorem~\ref{thm:recognize} gives a terminating recognition
algorithm. On a negative input the algorithm constructs generators of
an embedded \(\BS(1,n)\), \(n\geq2\), or \(\Z^2\). The corresponding
subgroup search is unconditional.

A finitely generated group is word-hyperbolic if geodesic triangles in
its Cayley graph are uniformly thin. Hyperbolic groups contain neither
\(\Z^2\) nor an infinite-order element conjugate to a proper power.
Lemma~\ref{lem:bs-types} therefore implies that they are BS-free.
Consequently, (H) identifies BS-freeness with hyperbolicity for
one-relator groups. This gives the isomorphism and conjugacy results
in Section~\ref{sec:consequences}, as well as decidability of equations.
For torsion-free groups in this class, decidability of the elementary
theory follows from the results on torsion-free hyperbolic groups
\cite{KMHyp,SelaVII}.

Section~\ref{sec:effective-vfc} gives an effective construction of a
virtual free-by-cyclic presentation for a hyperbolic one-relator group,
assuming the virtual compact specialness result stated as Theorem~1.1
in \cite{OAISpecial}. Theorem~\ref{vfc:thm:main} constructs a finite-index
normal subgroup $H$ and an embedding of $H$ in a mapping torus of a
finite-rank free-group automorphism. When the cyclic quotient is
nontrivial, it gives a finite HNN presentation of $H$ with finitely
generated free base and associated subgroups, an effective free basis
of the kernel, and the action of the stable letter on this basis.
The kernel may have infinite
rank. If a finite-rank free kernel is known to exist in a finite-index
subgroup, Proposition~\ref{vfc:prop:finite} gives a separate construction
by enumeration of finite presentations and Tietze transformations.

In Section~\ref{sec:toral-vfc}, we prove that, for a finitely generated
torsion-free virtually free-by-cyclic group, toral relative hyperbolicity
is equivalent to CSA, and also to the absence of $F_2\times\Z$ and the
Klein bottle group (Theorem~\ref{toral:thm:main}). Here the free kernel
need not be finitely generated. In the nonabelian case, the noncyclic
peripheral subgroups may be chosen isomorphic to $\Z^2$. In particular, every finitely generated
subgroup of a torsion-free CSA virtually free-by-cyclic group is toral
relatively hyperbolic. This result does not use either of the cited
results from \cite{OAI,OAISpecial}.

Section~\ref{sec:csa-algorithms} applies these results to decision
problems. Let $\mathcal C$ be a class of torsion-free one-relator groups
such that every CSA group in $\mathcal C$ is toral relatively hyperbolic.
Proposition~\ref{cases:prop:algorithm} gives a uniform algorithm on
$\mathcal C$ deciding whether a group is CSA and, if so, whether it
contains $\Z^2$. The algorithm runs a search for a failure of the CSA condition
in parallel with recognition of relative hyperbolicity with abelian
peripheral groups. We apply it to virtually free-by-cyclic groups,
ascending and partial ascending extensions of free groups, and finitely
generated subgroups of free-group mapping tori. Further applications to
graphs of groups and one-relator hierarchies assume the combination
statement of \cite{Stable}, with the hypotheses specified in
Section~\ref{sec:csa-algorithms}. These algorithms do not require a
procedure recognizing the class $\mathcal C$.

We shall use the theorem  about torsion in one-relator groups \cite[Chapter II]{LS}.
After free and cyclic reduction, a trivial relator gives a free group,
a nontrivial proper-power relator gives a group with torsion, and a
non-power relator gives a torsion-free group.

\section{First-order properties}\label{sec:fo}

\begin{definition}
A group is \emph{locally cyclic} if every finitely generated subgroup
is cyclic. A subgroup \(A\leq G\) is \emph{malnormal} if
\(A\cap gAg^{-1}=1\) for every \(g\notin A\).
A group is \emph{CSA} if its maximal abelian subgroups are malnormal.
Write \(C_G(x)=\{u:[x,u]=1\}\).
\end{definition}

We use the following results on subgroups of one-relator groups. An abelian subgroup of a
torsion-free one-relator group is either locally cyclic or isomorphic
to \(\Z^2\): this follows from
\cite[Lemma 5]{GKM}, since the only abelian group \(\BS(1,m)\)
is \(\BS(1,1)\). Also, a nonabelian torsion-free one-relator group
is CSA exactly when it contains neither \(F_2\times\Z\) nor a
nonabelian \(\BS(1,m)\), by \cite[Theorem 7]{GKM}.

\begin{lemma}\label{lem:bs-types}
Every Baumslag--Solitar group contains either \(\Z^2\) or
\(\BS(1,n)\) for an integer \(n\geq2\).
\end{lemma}
\begin{proof}
This is \cite[Lemma 2.10]{GKL}, with a sign adjustment.
If \(|m|,|n|>1\), the Bass--Serre normal form embeds
\(\langle a,b\mid a^n=b^m\rangle\), where \(b=tat^{-1}\).
Its central element \(z=a^n=b^m\) has infinite order, and its
quotient by \(\langle z\rangle\) is \(C_{|n|}*C_{|m|}\).
The image of \(ab\) has infinite order, so \(\langle z,ab\rangle\cong\Z^2\).
Otherwise the group is isomorphic to \(\BS(1,k)\).
For \(k=\pm1\) it contains \(\Z^2\); for \(k\geq2\) it is already
of the required type; for \(k\leq-2\), the subgroup generated by
the base element and the square of the stable letter is
\(\BS(1,k^2)\).
\end{proof}

Consider the following sentences in the  language of groups:
\begin{align}
\CT\;:=\;&\forall x,u,v\;
 \bigl(x\neq1\land[x,u]=[x,v]=1\ \Longrightarrow\ [u,v]=1\bigr),
 \label{eq:ct}\\
\Mal\;:=\;&\forall x,y,t\;
 \bigl(x\neq1\land y\neq1\land[x,y]=[x,tyt^{-1}]=1
 \ \Longrightarrow\ [x,t]=1\bigr),\label{eq:mal}\\
\Sq\;:=\;&\forall x,u,v\;
 \bigl(x\neq1\land[x,u]=[x,v]=1\ \Longrightarrow\nonumber\\
 &\hspace{2em}\exists w\;([x,w]=1\land
 (u=w^2\lor v=w^2\lor uv=w^2))\bigr).\label{eq:sq}
\end{align}
Here chaines of equalities abbreviate conjunctions. Put
\[
\CSA:=\CT\land\Mal,
\qquad \Phi:=\CSA\land\Sq.
\]

\begin{lemma}\label{lem:csa-formula}
For every group \(G\), \(G\models\CSA\) if and only if \(G\) is CSA.
If \(G\) is CSA, \(\Sq\) says that every nontrivial centralizer
\(A=C_G(x)\) satisfies \(\dim_{\mathbb F_2}(A/2A)\leq1\), in additive notation.
\end{lemma}
\begin{proof}
The sentence \(\CT\) says that every nontrivial centralizer is abelian;
these centralizers are then precisely the maximal abelian subgroups.
For \(A=C_G(x)\), a nontrivial element of \(A\cap tAt^{-1}\)
gives a \(y\neq1\) with \(y,tyt^{-1}\in A\). The sentence
\(\Mal\) forces \(t\in A\), and conversely malnormality gives this implication.

In an abelian group, \(\Sq\) states that, for any two classes
\(\bar u,\bar v\in A/2A\), one of \(\bar u,\bar v,\bar u+\bar v\)
is zero. This is precisely the absence of two independent vectors.
\end{proof}

\begin{theorem}\label{thm:fo}
For every torsion-free one-relator group \(G\),
\[
G\text{ is BS-free}\quad\Longleftrightarrow\quad G\models\Phi.
\]
\end{theorem}
\begin{proof}
Suppose first that \(G\) is BS-free. It contains neither
\(F_2\times\Z\), which contains \(\Z^2\), nor a nonabelian
\(\BS(1,m)\). Thus it is CSA by the structural criterion above;
the abelian cases are immediate. Each nontrivial centralizer is abelian
and cannot be \(\Z^2\), so is locally cyclic. For \(u,v\) in such
a centralizer, choose a generator \(c\) of \(\langle u,v\rangle\)
and write \(u=c^i,v=c^j\). One of \(i,j,i+j\) is even, giving the
required square root in the same centralizer. Hence \(G\models\Phi\).

Conversely, suppose \(G\models\Phi\). Its centralizers are abelian
and therefore either locally cyclic or \(\Z^2\). The latter possibility
violates \(\Sq\), by taking a basis of \(\Z^2\). Thus all nontrivial
centralizers are locally cyclic, and \(G\) contains no \(\Z^2\).
CSA groups contain no nonabelian \(\BS(1,n)\): its nontrivial abelian
base subgroup has a normalizer larger than its maximal abelian subgroup.
Lemma~\ref{lem:bs-types} now excludes every Baumslag--Solitar subgroup.
\end{proof}

Theorem~\ref{thm:fo} does not use either result from \cite{OAI,OAISpecial}.
Combining it with (H), the same sentence characterizes hyperbolicity
among torsion-free one-relator groups.

\begin{corollary}\label{cor:fo-partition}
Within torsion-free one-relator groups, the three sentences
\[
\CSA\land\Sq,\qquad \CSA\land\neg\Sq,\qquad \neg\CSA
\]
define a partition: respectively the BS-free groups, the CSA
groups containing \(\Z^2\), and the non-CSA groups.
BS-freeness is preserved under elementary equivalence within this class.
\end{corollary}

The recognition problem is considered in
Section~\ref{sec:recognition}. The first-order characterization alone
does not provide an algorithm for evaluating \(\Phi\) from a presentation.

\begin{corollary}\label{cor:all-fo}
The sentence
\[
\Psi:=\bigl(\exists z\;(z\neq1\land z^2=1)\bigr)\lor\Phi
\]
characterizes BS-freeness within all finitely generated one-relator groups.
Consequently, two one-relator groups, one BS-free and the other containing
a Baumslag--Solitar subgroup, cannot be elementarily equivalent.
\end{corollary}
\begin{proof}
For torsion-free inputs this is Theorem~\ref{thm:fo}.
A one-relator group with torsion is hyperbolic and hence BS-free.
If it has an involution, the first disjunct holds. Otherwise it contains
no infinite dihedral subgroup, so is CSA by \cite[Theorem 8]{GKM}.
Newman's centralizer theorem \cite[Theorem 2]{Newman} says that
its nontrivial centralizers are cyclic; hence \(\Sq\) holds too.
Elementary equivalence preserves the truth of \(\Psi\).
\end{proof}

\section{Recognition of BS-free one-relator groups}\label{sec:recognition}

We construct a recognition procedure using a Magnus hierarchy and
finite witnesses for Baumslag--Solitar subgroups. All tree actions are
simplicial and without inversions. An element is \emph{elliptic} if it
fixes a vertex and \emph{hyperbolic on the tree} if it translates along
an axis. Its translation length is
\(\length_T(g)=\min_{v\in V(T)}d_T(v,gv)\).

\begin{lemma}\label{lem:hierarchy}
For a nontrivial non-power relator, one can construct a finite hierarchy
\[
G=G_0\geq G_1\geq\cdots\geq G_N=1,
\qquad
G_i=\langle G_{i+1},t_i\mid t_iA_it_i^{-1}=B_i\rangle,
\]
where the groups \(G_i\) are explicitly embedded torsion-free one-relator
groups and \(A_i,B_i\) are free Magnus subgroups of \(G_{i+1}\).
The edge maps and all generator images in \(G\) are computable.
The word problem, membership in \(A_i,B_i\), and the test for
\(\length_{T_i}(g)>0\) are decidable.
\end{lemma}
\begin{proof}
Use Linton's construction of a one-relator hierarchy by \(\mathbb Z\)-domains in infinite cyclic covers
\cite{Linton}.
Its construction consists of finite graph operations in cyclic covers.
The cyclic character is obtained by integer linear algebra; domains,
overlaps, spanning trees, and their inclusion maps are computable.
One may enumerate admissible finite domains until the finite graph
conditions and the required strict decrease of complexity hold.
Let \(\Gamma\) be the underlying graph and
\(\lambda:S^1\to\Gamma\) the attaching loop, with combinatorial length
\(|\lambda|\). The decreasing nonnegative lexicographic complexity is
\[
\bigl(|\lambda|-|V(\lambda(S^1))|,\operatorname{rank}\pi_1\Gamma\bigr);
\]
see also the explicit construction in \cite[Lemma 8.3]{OAI}.
Existence at each nonterminal step and strict decrease ensure termination.
The Freiheitssatz and the HNN normal form show that the inclusions
in the hierarchy are injective. The terminal group is finite cyclic;
as it embeds in a torsion-free group, it is trivial.

Magnus's extended word-problem algorithm decides membership in the
edge groups \cite{LS}. When membership holds, an expression in their
generators can be found by enumeration and the word problem.
This makes Britton reduction and cyclic Britton reduction effective.
An element has positive tree translation length precisely when its
cyclically reduced HNN form contains a stable letter. 
\end{proof}

\begin{lemma}\label{lem:power}
In any group, if \(u\) has infinite order and
\(vuv^{-1}=u^n\), where \(|n|\geq2\), then
\(\langle u,v\rangle\cong\BS(1,n)\), with the indicated marking.
\end{lemma}
\begin{proof}
Write \(\BS(1,n)=\Z[1/n]\rtimes\Z\), with multiplication
\((q,k)(r,l)=(q+n^kr,k+l)\).
Every nontrivial normal subgroup meets the cyclic subgroup generated
by \((1,0)\): a nonzero base element has a nonzero integer multiple,
and the commutator of \((q,k)\), \(k\neq0\), with \((1,0)\) is
\((n^k-1,0)\neq(0,0)\).
The displayed relation gives a homomorphism to the ambient group.
A nontrivial kernel would kill a nonzero power of \(u\), which is impossible.
\end{proof}

\begin{lemma}\label{lem:tree}
At a level of the hierarchy, suppose
\[
a\in G_{i+1}\setminus\{1\},\qquad g\in G_i,
\qquad[a,g]=1,\qquad\length_{T_i}(g)>0.
\]
Then \(\langle a,g\rangle\cong\Z^2\).
\end{lemma}
\begin{proof}
If \(a^pg^q=1\), then \(g^q=a^{-p}\) is elliptic.
Since \(\length_{T_i}(g^q)=|q|\length_{T_i}(g)\), we get \(q=0\).
Torsion-freeness and \(a\neq1\) give \(p=0\).
\end{proof}

\begin{lemma}\label{lem:z2}
If \(G\) contains \(\Z^2\), a pair as in Lemma~\ref{lem:tree}
exists at some level of the computed hierarchy.
\end{lemma}
\begin{proof}
Choose the largest \(i\) such that \(G_i\) contains a subgroup
\(C\cong\Z^2\). The action of \(C\) on \(T_i\) cannot fix a vertex,
since vertex stabilizers are conjugates of \(G_{i+1}\).
A finitely generated group whose elements are all elliptic on a tree
has a fixed vertex. Thus \(C\) contains a translating element \(g\).
All of \(C\) preserves its axis, and no element reverses it, since
every element commutes with \(g\).
A translation defines a nonzero homomorphism \(C\to\Z\).
Its kernel has rank one. Choose \(a\neq1\) in the kernel;
it fixes the axis pointwise. Conjugating the pair puts an axis vertex
at the vertex stabilized by \(G_{i+1}\), so \(a\in G_{i+1}\).
\end{proof}

\begin{definition}\label{def:cert}
For a torsion-free input and its fixed computed hierarchy, a certificate
of type S is a tuple \((u,v,n)\), \(n\geq2\), satisfying
\[
u\neq1,\qquad vuv^{-1}=u^n.\tag{S}
\]
A certificate of type A is a tuple \((i,a,g)\) satisfying
\[
a\in G_{i+1}\setminus\{1\},\quad g\in G_i,
\quad[a,g]=1,\quad\length_{T_i}(g)>0.\tag{A}
\]
All words are written in the computed alphabets, with explicit images
in the original presentation.
\end{definition}

\begin{proposition}\label{prop:negative}
There is a semialgorithm that halts exactly on one-relator groups
containing a Baumslag--Solitar subgroup. When it halts, it gives
explicit generators of an embedded \(\BS(1,n)\), \(n\geq2\), or \(\Z^2\).
\end{proposition}
\begin{proof}
Trivial and proper-power relators define subgroups BS-free groups. On a nontrivial non-power input, construct the hierarchy and enumerate
all tuples in (S) and (A),  across the finitely many levels.
Every check terminates by Lemma~\ref{lem:hierarchy}.
It follows from Lemmas~\ref{lem:power} and~\ref{lem:tree} that the procedure halts only when \(G\) contains a Baumslag–Solitar subgroup.
Every Baumslag--Solitar subgroup contains one of the types in
Lemma~\ref{lem:bs-types}. An embedded solvable group supplies an
(S) tuple, and an embedded \(\Z^2\) supplies an (A) tuple by
Lemma~\ref{lem:z2}. Enumeration eventually reaches it.
\end{proof}

\begin{definition}\label{def:hyp}
For a finite presentation \(P\), \(\Hyp(P)\) is Papasoglu's
\emph{hyperbolicity-detection semialgorithm} \cite[pp.~193--200]{Papasoglu}.
On hyperbolic inputs it certifies hyperbolicity and computes a
thin-triangle constant \(\delta\); otherwise it need not terminate.
Dahmani--Guirardel \cite[Section 2.3.1]{DG} recall
the procedure.
\end{definition}

\begin{theorem}\label{thm:recognize}
Under (H), the following uniform procedure decides BS-freeness
and hyperbolicity for arbitrary finite one-relator presentations.
\begin{enumerate}
\item Freely and cyclically reduce \(r\). If \(r=1\), return
\texttt{BS-free}. If \(r\) is a nontrivial proper power, return
\texttt{BS-free}, by Newman.
\item For a remaining input, compute the hierarchy of
Lemma~\ref{lem:hierarchy}.
\item Run \(\Hyp(\langle X\mid r\rangle)\) and both certificate
enumerations (S) and (A) in parallel.
\item A hyperbolicity certificate gives the answer \texttt{BS-free}.
An (S) or (A) certificate gives the answer \texttt{contains BS},
together with the subgroup generators and their data that can be checked algorithmically.
\end{enumerate}
Hyperbolicity data can also be computed in every hyperbolic case.
\end{theorem}
\begin{proof}
Every answer returned by the procedure is correct, by Newman's theorem,
Lemmas~\ref{lem:power} and~\ref{lem:tree}, and the
fact that hyperbolic groups are BS-free. For a negative input,
Proposition~\ref{prop:negative} guarantees termination of the 
search. For a positive input, (H) guarantees that \(\Hyp\) terminates.
This proves termination on all inputs.
In a preliminary positive branch, one may additionally run \(\Hyp\)
to obtain \(\delta\).
\end{proof}

\begin{remark}
On an input containing a Baumslag--Solitar subgroup, the procedure
constructs an embedded \(\BS(1,n)\), \(n\geq2\), or \(\Z^2\).
Its parameters may differ from those of another Baumslag--Solitar subgroup
of the same group. The argument gives termination, without a bound on
the running time in terms of the relator length.
\end{remark}

\section{Algorithmic and model-theoretic consequences}\label{sec:consequences}

The corollaries below apply known theorems on hyperbolic groups to BS-free
one-relator groups using (H). This dependence concerns these applications;
the original model-theoretic results of the paper are independent of
both results cited from \cite{OAI,OAISpecial}.

\begin{corollary}\label{cor:iso}
Assume (H).
There is a uniform algorithm comparing a BS-free one-relator group
\(G\) with an arbitrary one-relator group \(H\), deciding whether
\(G\cong H\), and constructing an isomorphism when one exists.
\end{corollary}
\begin{proof}
Recognize \(H\). If it contains a BS subgroup, reject.
Otherwise both groups are hyperbolic, so apply
Dahmani--Guirardel's isomorphism algorithm \cite{DG}.

\end{proof}

\begin{corollary}\label{cor:equations}
Assume (H).
On the recognized class of BS-free one-relator groups, the conjugacy
and simultaneous conjugacy problems, and satisfiability of finite
systems of equations and inequations with coefficients, are uniformly
decidable. These statements include the torsion case.
\end{corollary}
\begin{proof}
Compute hyperbolicity data and apply the hyperbolic-group algorithms
in \cite[Section 2.3]{DG} and \cite{DGEquations}.
\end{proof}

\begin{corollary}\label{cor:model}
Assume (H).
Every torsion-free BS-free one-relator group has decidable elementary
theory and effective reduction of formulas to Boolean combinations
of \(\forall\exists\)-formulas. It is equationally Noetherian.
Its nontrivial centralizers are infinite cyclic, it is automatic,
and it is Hopfian.
\end{corollary}
\begin{proof}
Follows the corresponding torsion-free hyperbolic-group results:
\cite{KMHyp,SelaVII} for elementary theory and algebraic geometry,
the standard hyperbolic centralizer and automaticity theorems,
and \cite{SelaHopf} for Hopficity.
\end{proof}

A group is \emph{equationally Noetherian} if every system of equations
in finitely many variables has the same solution set as a finite
subsystem. Under (H), the elementary-equivalence results for torsion-free
hyperbolic groups \cite{SelaVII} apply to torsion-free BS-free one-relator
groups. By Corollary~\ref{cor:all-fo}, every one-relator group elementarily
equivalent to a BS-free one-relator group is itself BS-free.

\begin{corollary}\label{cor:mutual}
Assume (H).
If \(G\) is a torsion-free BS-free one-relator group and there are
epimorphisms \(G\to H\) and \(H\to G\), then \(G\cong H\).
\end{corollary}
\begin{proof}
The composite is a surjective endomorphism of the Hopfian group \(G\)
and hence injective. Thus the epimorphism \(G\to H\) is also injective.
\end{proof}

\begin{remark}\label{rem:special}
If a hyperbolic group has a finite quasiconvex hierarchy in the sense
of Wise, then it is virtually compact special
\cite{Wise}. It is therefore residually
finite and has a finite-index subgroup which embeds in a right-angled
Artin group. 
Here the edge groups of the hierarchy must be quasiconvex in the
ambient group at each splitting.
\end{remark}

\section{Interpretation of solvable Baumslag--Solitar groups}
\label{sec:interpretation}

We recall the definition of interpretability.

\begin{definition}
An algebraic structure \(B\) is \emph{interpretable} in a structure \(G\)
if there is a definable set \(\mathcal D\subseteq G^r\), a definable
equivalence relation \(E\) on \(\mathcal D\), and definable operations and
relations on \(\mathcal D/E\) which make this quotient isomorphic to \(B\).
The formulas may contain a finite tuple of parameters.
The interpretation is \emph{regular} if these parameters range over a
nonempty definable set \(P\), and the same formulas give a copy of \(B\)
for every tuple in \(P\). It is \emph{injective} if \(E\) is equality.
An interpretation in a class of structures is \emph{uniform} if the
defining formulas are the same for every structure in the class. See \cite[Section 4.2]{KMS}.
\end{definition}

\begin{lemma}\label{lem:subgroup-trichotomy}
Every subgroup of a torsion-free one-relator group is locally cyclic,
contains a nonabelian free group, or is isomorphic to \(\BS(1,n)\) for a
nonzero integer \(n\). In particular, every abelian subgroup is locally
cyclic or isomorphic to \(\Z^2\).
\end{lemma}
\begin{proof}
This is \cite[Lemma 5]{GKM}; no finite-generation assumption on the
subgroup is required.
\end{proof}

\begin{lemma}
\label{lem:loccyc-centralizer}
Let \(G\) be a torsion-free one-relator group and let \(A\leq G\) be
noncyclic and locally cyclic. Then \(C_G(A)\) is locally cyclic.
\end{lemma}
\begin{proof}
Put \(L=C_G(A)\). We have \(A\leq Z(L)\). If \(L\) contained a free
group \(F\) of rank two, then \(A\cap F=1\), since this intersection
is contained in \(Z(F)\). For a nontrivial \(f\in F\), the subgroup
\(A\langle f\rangle\) would therefore be \(A\times\Z\). This is an
abelian subgroup which is neither locally cyclic nor \(\Z^2\):
it has rank two and is not finitely generated. This contradicts
the abelian-subgroup classification.

By \cite[Lemma 5]{GKM}, \(L\) is locally cyclic or
\(\BS(1,n)\). The latter is impossible. For \(|n|>1\) its center
is trivial; for \(n=-1\) its center is infinite cyclic; and for \(n=1\)
it is \(\Z^2\). None of these centers can contain the noncyclic locally
cyclic group \(A\). Thus \(L\) is locally cyclic.
\end{proof}

\begin{lemma}
\label{lem:definable-base}
Let \(G\) be a torsion-free one-relator group, let \(|m|\geq2\), and suppose
\[
a\neq1,\qquad tat^{-1}=a^m.
\]
Define
\[
D=D_m(a,t)=\{x\in G:[x,a]=1,\;txt^{-1}=x^m\}.
\]
Then \(D\) is a noncyclic locally cyclic subgroup, \(tDt^{-1}=D\),
and \(C_G(D)=D\).
\end{lemma}
\begin{proof}
By Lemma~\ref{lem:power}, the pair \(a,t\) generates a copy of
\(\BS(1,m)\). Its abelian base is
\[
A=\bigcup_{j\geq0}t^{-j}\langle a\rangle t^j
  \cong\Z[1/|m|].
\]
The cyclic groups form an increasing chain; \(A\) is noncyclic and
normalized by \(t\). Lemma~\ref{lem:loccyc-centralizer} gives that
\(L=C_G(A)\) is locally cyclic, and \(tLt^{-1}=L\).

We claim \(D=L\). If \(x\in D\), then for every \(j\geq0\),
\[
t^j[x,t^{-j}at^j]t^{-j}=[x^{m^j},a]=1.
\]
Thus \(x\) centralizes all of \(A\), so \(D\subseteq L\).
Conversely, if \(x\in L\), choose integers \(p\neq0,q\) with
\(x^p=a^q\), using local cyclicity. Both \(txt^{-1}\) and \(x^m\)
belong to \(L\), and
\[
(txt^{-1})^p=a^{mq}=(x^m)^p.
\]
Unique roots in the torsion-free abelian group \(L\) imply
\(txt^{-1}=x^m\). Hence \(x\in D\).

Since \(A\leq D\) and \(D\) is abelian,
\[
D\leq C_G(D)\leq C_G(A)=D,
\]
which proves self-centralization. The other assertions follow from
\(D=L\).
\end{proof}

\begin{proposition}\label{prop:definable-normalizer}
With the notation above, put \(N=N_G(D)\). Then for some integer
\(n\) with \(|n|\geq2\),
\[
N\cong\BS(1,n),\qquad
D\cong\Z[1/n],\qquad N/D\cong\Z.
\]
Under the standard orientation of this quotient, the image of \(t\)
is an integer \(k\geq1\), and
\[
n^k=m,\qquad k\leq\lfloor\log_2|m|\rfloor.
\]
Moreover, \(\langle D,t\rangle\cong\BS(1,m)\) and has index \(k\) in \(N\).
\end{proposition}
\begin{proof}
Conjugation gives \(N/D\hookrightarrow\operatorname{Aut}(D)\), because
\(C_G(D)=D\). A torsion-free locally cyclic group embeds in \((\Q,+)\);
its automorphisms are multiplications by nonzero rationals. Its
automorphism group is therefore abelian. Hence \(N\) is metabelian.
It contains the nonabelian group \(\langle a,t\rangle\), so \cite[Lemma 5]{GKM} gives \(N\cong\BS(1,n)\). The cases \(n=1,-1\) are impossible:
these groups are virtually \(\Z^2\), all their subgroups are finitely
generated, and \(D\) is noncyclic locally cyclic.

In \(\mathbb Z[1/n]\rtimes\mathbb Z\), every element outside the normal abelian subgroup \(\mathbb Z[1/n]\times\{0\}\) has cyclic centralizer. The projection of its centralizer to \(\Z\) is injective,
since multiplication by \(n^h\), \(h\neq0\), fixes no nonzero base element.
Thus the noncyclic abelian subgroup \(D\) lies in the standard base \(B\).
Conversely, \(B\) centralizes \(D\), so \(B\leq C_N(D)=D\).
Consequently \(B=D\).

If \(t\) has height \(k\), its conjugation action on \(D\) is multiplication
by \(n^k\). Lemma~\ref{lem:definable-base} identifies this action as
multiplication by \(m\), so \(n^k=m\). Since \(|n|\geq2\) and \(|m|\geq2\),
we have \(k\geq1\) and \(k\leq\lfloor\log_2|m|\rfloor\).
Finally \(\Z[1/n]=\Z[1/m]\), and
\(\langle D,t\rangle=D\rtimes\langle t\rangle\), with action \(m\).
This is \(\BS(1,m)\); its quotient image in \(N/D\) is \(k\Z\).
\end{proof}

\begin{definition}
\label{def:interpretation-formulas}
Fix \(|m|\geq2\) and put \(K=\lfloor\log_2|m|\rfloor\). All exponents below
are fixed integers, so these are ordinary first-order group formulas.
Define the parameter condition and the base and normalizer formulas by
\begin{align*}
P_m(a,t)&:=(a\neq1\ \land\ tat^{-1}=a^m),\\
d_m(x;a,t)&:=([x,a]=1\ \land\ txt^{-1}=x^m),\\
\nu_m(g;a,t)&:=\forall x\,
 \bigl(d_m(x;a,t)\ \longleftrightarrow\
       d_m(gxg^{-1};a,t)\bigr).
\end{align*}
For \(1\leq j\leq K\), put
\[
R_j(z;a,t):=\nu_m(z;a,t)\ \land\
 \exists y\,\bigl(\nu_m(y;a,t)\land d_m(zy^{-j};a,t)\bigr).
\]
For \(1\leq k\leq K\), put
\[
\eta_k(a,t):=R_k(t;a,t)\ \land\
 \bigwedge_{j=k+1}^{K}\neg R_j(t;a,t),
\]
where an empty conjunction is true. Finally, define the group domain
\[
\Delta_m(g;a,t):=
 \bigvee_{k=1}^{K}\bigl(\eta_k(a,t)\land R_k(g;a,t)\bigr).
\]
\end{definition}

\begin{theorem}\label{thm:bs-interpretation}
Let \(m\) be an integer with \(|m|>1\), and let \(G\) be a one-relator
group containing \(\BS(1,m)\). Then \(\BS(1,m)\) is regularly and
injectively interpretable in \(G\). The interpretation is uniform for
fixed \(m\), with admissible parameter set \(P_m(a,t)\).
For every admissible pair, the formula \(\Delta_m(g;a,t)\) defines
a subgroup isomorphic to \(\BS(1,m)\), whose operations are restrictions
of the operations of \(G\).
The interpretation formulas can be effectively constructed from \(m\).
\end{theorem}
\begin{proof}
The group \(G\) is torsion-free, since one-relator groups with torsion
are hyperbolic and contain no Baumslag--Solitar subgroup. An embedded
\(\BS(1,m)\) supplies a pair satisfying \(P_m\). Every pair satisfying
\(P_m\) generates such a subgroup by
Lemma~\ref{lem:power}. Fix any such pair. Lemma~\ref{lem:definable-base}
and Proposition~\ref{prop:definable-normalizer} apply to it.
The formulas \(d_m,\nu_m\) define \(D,N\), respectively.

Write \(h:N\to N/D\cong\Z\) for the height map, oriented as in
Proposition~\ref{prop:definable-normalizer}, and let \(h(t)=k\).
For \(z\in N\),
\[
R_j(z;a,t)\quad\Longleftrightarrow\quad h(z)\in j\Z.
\]
Indeed \(zy^{-j}\in D\) is precisely \(h(z)=jh(y)\), and every integer
is a height of an element of \(N\). Thus \(R_j(t;a,t)\) holds exactly
when \(j\mid k\). Since \(1\leq k\leq K\), the largest \(j\leq K\)
for which this holds is \(k\) itself. Exactly one formula \(\eta_j\)
holds, namely \(\eta_k\).

It follows that \(\Delta_m(G;a,t)=h^{-1}(k\Z)=\langle D,t\rangle\),
which is isomorphic to \(\BS(1,m)\). This works for every admissible
pair. Equality is ordinary equality and all group operations are
restrictions of the ambient operations; hence this is an injective
interpretation. The finite formulas can plainly be constructed from \(m\).
\end{proof}

\begin{corollary}
\label{cor:containing-undecidability}
If a one-relator group \(G\) contains \(\BS(1,m)\) with \(|m|>1\),
then arithmetic is regularly interpretable in \(G\), and its elementary
theory in the language of groups is undecidable.
\end{corollary}
\begin{proof}
Such a group is torsion-free: a one-relator group with torsion is
hyperbolic and contains no Baumslag--Solitar subgroup. If \(m<-1\),
replace the stable letter by its square to obtain \(\BS(1,m^2)\).
It is therefore enough to treat a fixed \(m\geq2\).

By Theorem~\ref{thm:bs-interpretation}, \(\BS(1,m)\) is regularly
interpretable in \(G\). The regular interpretation of the ring \(\Z\)
in \(\BS(1,m)\), given in \cite[Theorem 3 and Section 4.4]{DMBS},
composes with this interpretation. Its parameter conditions and
quantifiers are relativized to \(\Delta_m\). The resulting admissible
family is first-order definable and nonempty.

Let \(\sigma\) be a sentence in the language of groups. Relativizing
its quantifiers to \(\Delta_m\) gives a formula
\(\sigma^{\Delta_m}(a,t)\). Then
\[
\BS(1,m)\models\sigma
\quad\Longleftrightarrow\quad
G\models\exists a\,\exists t\,
 \bigl(P_m(a,t)\land\sigma^{\Delta_m}(a,t)\bigr).
\]
Every admissible pair gives the same abstract group, so this equivalence
is valid even for sentences with arbitrary quantifier alternation.
It reduces the undecidable theory of \(\BS(1,m)\) to \(\Th(G)\).
\end{proof}

\begin{example}
Let
\[
G=\langle b,t\mid tb^2t^{-1}=b^4\rangle,\qquad a=b^2.
\]
The elements \(b\) and \(c=t^{-1}b^2t\) belong to \(C_G(a)\), and
Britton's lemma gives \([b,c]\neq1\). Thus \(C_G(a)\) is nonabelian.
On the other hand, \(c^2=a\) and \(tct^{-1}=c^2\), so
\(c\in D_2(a,t)\). By Lemma~\ref{lem:definable-base}, \(D_2(a,t)\)
is locally cyclic and normalized by \(t\). Since \(K=1\), the formula
\(\Delta_2\) defines \(N_G(D_2(a,t))\), which is isomorphic to
\(\BS(1,2)\).
\end{example}

\begin{example}
In \(G=\langle b,t\mid tbt^{-1}=b^2\rangle\), put \(a=b^3\).
Then \(\langle a,t\rangle\) has index three in \(G\), whereas
\(D_2(a,t)\) is the full abelian base and \(\Delta_2\) defines \(G\).
Thus the definable subgroup given by
Theorem~\ref{thm:bs-interpretation} may properly contain the initially
given copy of \(\BS(1,m)\).
\end{example}

\begin{remark}
The groups \(\BS(1,1)=\Z^2\) and \(\BS(1,-1)\) have decidable elementary
theories. For the Klein bottle group \(\BS(1,-1)\), this follows by
expressing multiplication in its two integer coordinates in Presburger
arithmetic. The condition \(|m|>1\) in
Corollary~\ref{cor:containing-undecidability} is therefore necessary.
\end{remark}

\section{Uniform arithmetic in groups of cohomological dimension two}
\label{sec:nonuni-bs}

All cohomological dimensions in this section are over $\Z$.
Two subgroups are \emph{commensurable} if their intersection has finite
index in both. For an infinite cyclic subgroup $D\leq U$, put
\[
 \Comm_U(D)=\{g\in U:D\text{ and }gDg^{-1}
                         \text{ are commensurable}\}.
\]
This is a subgroup. We say that $D$ is \emph{commensurated by $U$} if
$\Comm_U(D)=U$. For infinite cyclic groups, commensurability is equivalent
to nontrivial intersection.

\begin{definition}
A generalized Baumslag--Solitar group, or GBS group, is the fundamental
group of a finite connected graph of infinite cyclic groups with
injective edge maps.
\end{definition}

\begin{theorem}\label{arith:thm:criterion}
Let $\mathcal C$ be the class of countable groups $U$ with
$\cdim U\leq2$ for which there are $a,s\in U$ and nonzero integers $m,n$
such that $\langle a\rangle$ is infinite cyclic and commensurated by $U$,
and
\[
 sa^ms^{-1}=a^n,\qquad |m|\ne|n|.
\]
The ring $(\Z,+,\times,0,1)$ is uniformly, regularly and injectively
interpretable in $\mathcal C$. In particular, the interpretation formulas
and the formula defining the admissible parameters are independent of
$m,n$ and of $U$.
\end{theorem}

We prove the theorem by defining $(\Z,+,\mid,0,1)$. Multiplication is
first-order definable in this structure by the classical result of
Julia Robinson \cite{JR}. We use the usual convention that every integer
divides zero and that zero divides only zero.

\subsection{Subgroups and centralizers}

We shall use the following two structural results.

\begin{theorem}\label{ambient:thm:kropholler}
A finitely generated group of cohomological dimension two with a
commensurated infinite cyclic subgroup is a GBS group.
\end{theorem}

This is the form of \cite[Theorem C]{Kropholler} stated in
\cite[proof of Theorem 3.2]{BK}.

\begin{theorem}\label{arith:thm:soluble}
A countable soluble group of cohomological dimension at most two is
locally cyclic or isomorphic to $\BS(1,k)$ for a nonzero integer $k$.
In particular, a countable abelian group of cohomological dimension at
most two is locally cyclic or isomorphic to $\Z^2$.
\end{theorem}

For dimension two this is \cite[Theorem 5]{Gildenhuys}; see also
\cite[Theorem 3 and its proof]{KLL}. Dimension one follows from the
Stallings--Swan theorem \cite[Chapter VIII]{Brown}.

\begin{lemma}\label{arith:lem:modular}
Suppose $\langle a\rangle$ is infinite cyclic and commensurated by $U$.
There is a homomorphism
\[
 \Delta:U\longrightarrow\Q^\times,\qquad
 ga^rg^{-1}=a^q\ \Longrightarrow\ \Delta(g)=q/r
 \quad(r,q\ne0).
\]
It is unchanged when $\langle a\rangle$ is replaced by a commensurable
infinite cyclic subgroup. Moreover,
\[
 \ker\Delta=\bigcup_{r\geq1}C_U(a^r),
\]
where the union is directed by divisibility of $r$.
\end{lemma}
\begin{proof}
Two conjugation relations give the same ratio after taking common
powers. Applying two conjugations to a sufficiently large common power
proves multiplicativity. Comparing common powers also proves the
assertion about commensurable cyclic subgroups. Finally, $\Delta(g)=1$
if and only if $g$ centralizes a nonzero power of $a$, and
$C_U(a^r),C_U(a^q)\leq C_U(a^{rq})$ for $r,q>0$.
\end{proof}

The homomorphism $\Delta$ is the \emph{modular homomorphism}.
In a GBS tree, adjacent vertex stabilizers contain their common edge
stabilizer with finite index. Following a finite path shows that all
vertex stabilizers are commensurable. They are therefore commensurated
by the group and define the same modular homomorphism. A GBS group is
\emph{non-unimodular} if $\Delta(U)\not\subseteq\{1,-1\}$; see
\cite[Section 2]{Levitt}.

If $sa^ms^{-1}=a^n$ with $|m|\ne|n|$, then $a$ is elliptic in every
action without inversions on a simplicial tree: translation length
gives $|m|\ell(a)=|n|\ell(a)$. In a GBS tree, every nontrivial elliptic
cyclic subgroup is consequently commensurable with $\langle a\rangle$.
Every elliptic element has modular value one, since its cyclic vertex
stabilizer contains and centralizes a nonzero power of $a$.

\begin{lemma}\label{arith:lem:local}
Let $U,a,s,m,n$ satisfy the hypotheses of
Theorem~\ref{arith:thm:criterion}. Every finitely generated subgroup
containing $a,s$ is a non-unimodular GBS group. Every subgroup of $U$
is locally cyclic, contains a nonabelian free group, or is isomorphic
to $\BS(1,k)$.
\end{lemma}
\begin{proof}
Let $J$ be finitely generated and contain $a,s$. It commensurates
$\langle a\rangle$ and has cohomological dimension at most two. It is
not free, since unequal powers of a nontrivial element of a free group
cannot be conjugate, by translation length in its Cayley tree.
Thus $\cdim J=2$, and Theorem~\ref{ambient:thm:kropholler} makes $J$
a GBS group. The preceding observation makes $a$ elliptic in this tree.
Its modular homomorphism is the restriction of $\Delta$, and
$|\Delta(s)|\ne1$.

Let $K\leq U$ contain no nonabelian free subgroup. Every finitely
generated $K_0\leq K$ lies in such a group $J$. The tree-action
alternative \cite{Serre} implies that $K_0$ fixes a vertex, fixes an end,
or preserves a line. In the first case $K_0$ is cyclic. In the second,
the kernel of the homomorphism to $\Z$ given by displacement along a
ray is locally cyclic: finitely many elements in the kernel fix a common
tail and lie in a cyclic vertex stabilizer. In the third case the
pointwise stabilizer of the line is cyclic and the image is a subgroup
of the infinite dihedral group. Thus $K_0$ is soluble in each case.
Theorem~\ref{arith:thm:soluble} makes it metabelian. Since metabelianity
is a group identity, $K$ is metabelian. Apply
Theorem~\ref{arith:thm:soluble} to $K$.
\end{proof}

\begin{lemma}\label{ambient:lem:center}
Let $V$ be either a torsion-free one-relator group or a group in
$\mathcal C$. If $z\ne1$ and $C_V(z)$ is nonabelian, then
$D=Z(C_V(z))$ is infinite cyclic. It is defined by
\begin{equation}\label{ambient:eq:delta}
 \delta(x;z):=[x,z]=1\ \land\
 \forall y\,([y,z]=1\ \Longrightarrow\ [x,y]=1).
\end{equation}
\end{lemma}
\begin{proof}
Put $C=C_V(z)$, so $1\ne z\in D$. If $C$ contains a nonabelian free
subgroup $F$, then $D\cap F\leq Z(F)=1$. For $1\ne f\in F$, the
abelian subgroup $D\langle f\rangle=D\times\langle f\rangle$ contains
$\langle z,f\rangle\cong\Z^2$. Theorem~\ref{arith:thm:soluble} identifies
it with $\Z^2$, so its direct factor $D$ is infinite cyclic.
Otherwise Lemma~\ref{lem:subgroup-trichotomy} or
Lemma~\ref{arith:lem:local}, respectively, gives $C\cong\BS(1,k)$.
Among the nonabelian groups of this form only $\BS(1,-1)$ has nontrivial
center, and its center is infinite cyclic. The  formula above is
the definition of $Z(C_V(z))$.
\end{proof}

\begin{lemma}\label{arith:lem:find-center}
If $U\in\mathcal C$ is not soluble, some $z=a^r$, $r\geq1$, has
nonabelian centralizer. The group $Z(C_U(z))$ is commensurable with
$\langle a\rangle$ and commensurated by $U$.
\end{lemma}
\begin{proof}
If every $C_U(a^r)$ were abelian, the directed union
$\ker\Delta=\bigcup_{r\geq1}C_U(a^r)$ would be abelian. Since
$U/\ker\Delta$ is abelian, $U$ would be metabelian.
For the chosen $z$, the cyclic group $Z(C_U(z))$ contains $a^r$.
This proves the remaining assertions.
\end{proof}

\begin{lemma}\label{ambient:lem:cyclic-T}
Let $U\in\mathcal C$. If $|\Delta(s)|\ne1$, then $C_U(s)$ is infinite
cyclic.
\end{lemma}
\begin{proof}
For $x\in C_U(s)$, put $J=\langle a,s,x\rangle$. The unequal-power
relation supplied by $\Delta(s)$ and Lemma~\ref{arith:lem:local} show
that $J$ is GBS. The element $a$ is elliptic and $s$ is hyperbolic.
The centralizer $C_J(s)$ preserves the axis of $s$ and acts by
translations. This action is faithful. Indeed, a nontrivial element $w$
in its kernel would be elliptic, so $w^p=a^q$ for nonzero $p,q$.
As $w$ commutes with $s$, this would give $\Delta(s)=1$.
Thus $C_J(s)$ embeds in the integral translation group of the axis
and is infinite cyclic.

Choose a prime $\ell$ with $v_\ell(\Delta(s))\ne0$. The homomorphism
\[
 v_\ell\circ\Delta:C_U(s)\longrightarrow\Z
\]
is injective: for each $x\in C_U(s)$ its restriction to the cyclic
group $C_J(s)$ is nonzero on $s$, and hence injective. Its image is a
nonzero subgroup of $\Z$, proving the lemma.
\end{proof}

\subsection{Intersections with conjugates of cyclic subgroups}

\begin{lemma}\label{arith:lem:cyclic-intersections}
Let $J$ be a GBS group, let $D=\langle d\rangle$ be a nontrivial
elliptic cyclic subgroup, and suppose
\[
 \Delta(r)=\varepsilon q/p,\qquad
 \varepsilon\in\{1,-1\},\quad p,q\geq1,\quad
 \gcd(p,q)=1,\quad p\ne q.
\]
There is an integer $e\geq1$ such that for every $n\geq1$,
\begin{align}
 D\cap r^{-n}Dr^n&=\langle d^{ep^n}\rangle,&
 D\cap r^nDr^{-n}&=\langle d^{eq^n}\rangle,
 \label{arith:eq:intersections}\\
 r^nd^{ep^n}r^{-n}&=d^{\varepsilon^n eq^n}.&&
 \label{arith:eq:conjugation}
\end{align}
In particular, for every nonzero integer $n$,
\begin{equation}\label{arith:eq:product-independent}
 (D\cap r^{-n}Dr^n)(D\cap r^nDr^{-n})=\langle d^e\rangle.
\end{equation}
\end{lemma}
\begin{proof}
The element $r$ is hyperbolic because its modular value has absolute
value different from one. Let $v$ be a vertex fixed by $D$, and let
$x$ be its projection onto the axis of $r$. Put
\[
 V=\operatorname{Stab}_J(x)=\langle b\rangle,\qquad
 D\cap V=\langle d^h\rangle=\langle b^\ell\rangle,
\]
where $h,\ell>0$, after choosing the orientation of $b$.
Commensurability of vertex stabilizers makes this intersection
nontrivial. For some $g\geq1$,
\[
 V\cap r^{-1}Vr=\langle b^{gp}\rangle,\qquad
 rb^{gp}r^{-1}=b^{\varepsilon gq}.
\]
The image is $V\cap rVr^{-1}=\langle b^{gq}\rangle$.

The vertices $x,r^{-1}x,\ldots,r^{-n}x$ lie in this order on the axis.
An element fixing the first and last fixes all the intermediate
vertices. Therefore, if $b^j$ and $r^nb^jr^{-n}$ belong to $V$, each
intermediate conjugate belongs to $V$. Iterating the displayed
one-step conjugation, and using $\gcd(p,q)=1$, gives
\[
 V\cap r^{-n}Vr^n=\langle b^{gp^n}\rangle,\qquad
 r^nb^{gp^n}r^{-n}=b^{\varepsilon^n gq^n}.
\]

Suppose $d^j\in D$ and $r^nd^jr^{-n}\in D$. The first element fixes
$v$ and $r^{-n}v$, hence fixes $x$; the conjugate fixes $v$ and $r^nv$,
hence also fixes $x$. Both endpoints therefore belong to $D\cap V$.
For an input $d^{hk}=b^{\ell k}$, the necessary and sufficient
conditions are
\[
 gp^n\mid\ell k,\qquad p^n\mid k.
\]
The first condition is membership in the domain inside $V$. The
second says that the output exponent
$\varepsilon^n\ell kq^n/p^n$ is a multiple of $\ell$, so that the
output lies in $D\cap V$.
Writing $k=p^nu$ reduces the first condition to $g\mid\ell u$.
Thus \eqref{arith:eq:intersections} and
\eqref{arith:eq:conjugation} hold with
\[
 e=h\,g/\gcd(g,\ell).
\]
Finally, $\gcd(p^n,q^n)=1$ gives
\eqref{arith:eq:product-independent} for positive $n$; negative $n$
interchange the two factors.
\end{proof}

\subsection{Divisibility of exponents}

\begin{lemma}\label{arith:lem:gcd}
Let $q>p\geq1$ be coprime integers and put $F_j=q^j-p^j$ for $j\geq0$.
For positive $i,j$,
\[
 \gcd(F_i,F_j)=F_{\gcd(i,j)}.
\]
For arbitrary integers $i,j$,
\[
 i\mid j\quad\Longleftrightarrow\quad F_{|i|}\mid F_{|j|}.
\]
\end{lemma}
\begin{proof}
For $j\geq i\geq1$,
\[
 F_j-q^{j-i}F_i=p^iF_{j-i}.
\]
Since $\gcd(F_i,p)=1$, the Euclidean algorithm on the indices proves
the gcd formula. The sequence $F_j$ is strictly increasing, with
$F_0=0$. Comparing a term with its gcd with another term gives the
divisibility assertion for positive indices. Absolute values and
$F_0=0$ give the remaining cases.
\end{proof}

\subsection{The nonsoluble case}

For parameters $z,s\in U$, let
\[
 D=\delta(U;z),\qquad
 R_s=\{x\in D:sxs^{-1}\in D\},\qquad
 S_s=\{x\in D:s^{-1}xs\in D\}.
\]
Thus $R_s=D\cap s^{-1}Ds$ and $S_s=D\cap sDs^{-1}$.
Let $\mathcal P_{\rm n}(z,s)$ be the conjunction of the following
first-order conditions:
\begin{enumerate}
\item $z\ne1$ and $C_U(z)$ is nonabelian;
\item $\forall g\,\exists x\,
 (x\ne1\land\delta(x;z)\land\delta(gxg^{-1};z))$;
\item $R_s\ne S_s$.
\end{enumerate}
The third condition is inequality of two definable sets, so it is
also first-order. By Lemma~\ref{ambient:lem:center}, $D$ is infinite
cyclic. The second condition says that $D$ is commensurated by $U$.
Write $\Delta_D$ for the modular homomorphism defined using $D$.
The third condition says that $|\Delta_D(s)|\ne1$, since conjugation by $s$
maps $R_s$ isomorphically onto $S_s$.

Assume $\mathcal P_{\rm n}(z,s)$ and define
\begin{equation}\label{arith:eq:nonsoluble-A}
 A=R_sS_s,\qquad T=C_U(s),\qquad
 E=T^2=\{v^2:[v,s]=1\}.
\end{equation}
The product defining $A$ is a subgroup because both factors lie in the
cyclic group $D$. All three sets are definable by fixed formulas.

For verification, write $D=\langle d\rangle$ and
$T=\langle r\rangle$, using Lemma~\ref{ambient:lem:cyclic-T} with base
subgroup $D$. Orient $r$ so that
\[
 \Delta_D(r)=\varepsilon q/p,\qquad q>p\geq1,\quad\gcd(p,q)=1.
\]
Since $s=r^k$ for some $k\ne0$, the finitely generated group
$\langle d,r\rangle$ is GBS and contains $s$.
Lemma~\ref{arith:lem:cyclic-intersections} gives
\[
 A=\langle d^e\rangle.
\]
Put $c=d^e$. For $n\geq0$,
\[
 r^nc^{p^n}r^{-n}=c^{\varepsilon^nq^n},\qquad
 r^nc^jr^{-n}\in A\ \Longleftrightarrow\ p^n\mid j.
\]
These integers and generators are used only in the proof; they do not
occur in the defining formulas.

For $y\in E$, put
\begin{equation}\label{arith:eq:I}
 I_y=\{yxy^{-1}x^{-1}:x\in A,\ yxy^{-1}\in A\}.
\end{equation}
If $\alpha(x;z,s)$ defines $A$, membership in $I_y$ is defined by
\begin{equation}\label{arith:eq:iota}
 \iota(w,y;z,s):=
 \exists x\,\bigl(\alpha(x;z,s)\land
 \alpha(yxy^{-1};z,s)\land w=yxy^{-1}x^{-1}\bigr).
\end{equation}
Since $E=\langle r^2\rangle$, the preceding calculation gives
\begin{equation}\label{arith:eq:I-values}
 I_{r^{2j}}=\langle c^{q^{2|j|}-p^{2|j|}}\rangle
 \quad(j\in\Z).
\end{equation}
For $j=0$ this means $I_1=\{1\}$.
Define, on $E$,
\begin{equation}\label{arith:eq:Div}
 \Div(y_1,y_2;z,s)\quad\Longleftrightarrow\quad
 I_{y_2}\subseteq I_{y_1}.
\end{equation}
Subgroup inclusion is expressed by
$\forall w\,(\iota(w,y_2;z,s)\Rightarrow\iota(w,y_1;z,s))$.
Lemma~\ref{arith:lem:gcd}, applied to $p^2,q^2$, proves
\begin{equation}\label{arith:eq:integer-divisibility}
 \Div(r^{2i},r^{2j};z,s)\quad\Longleftrightarrow\quad i\mid j.
\end{equation}

\subsection{The soluble case}

If $U\in\mathcal C$ is soluble, Theorem~\ref{arith:thm:soluble} gives
\[
 U\cong\BS(1,k)=B\rtimes\langle v\rangle,\qquad
 B=\Z[1/|k|],\quad |k|>1.
\]
Here conjugation by $v$ multiplies the additive group $B$ by $k$.
The projection $U\to U/B\cong\Z$ is the \emph{height map};
the height of $bv^j$ is $j$.
Locally cyclic groups and $\Z^2$ have no unequal-power relation.
Neither does the Klein bottle group: its normal subgroup $\Z^2$ of
index two has conjugation actions of order at most two, so such a
relation would imply $m^2=n^2$.

Let $\mathcal P_{\rm s}(z,s)$ say
\[
 z\ne1,\qquad C_U(z)\text{ is abelian and normal in }U,\qquad [z,s]\ne1.
\]
Normality here is the first-order condition
\[
 \forall g\,\forall x\,
 ([x,z]=1\ \Longleftrightarrow\ [gxg^{-1},z]=1).
\]
These conditions force $C_U(z)=B$ and $s\notin B$. Indeed, an element of
nonzero height has centralizer intersecting $B$ trivially, since its
conjugation acts on $B$ by a scalar different from one. If that
centralizer were normal, its commutators with $B$ would lie in this
trivial intersection, a contradiction. Thus $z\in B\setminus\{1\}$,
and its centralizer is exactly $B$.

Set
\[
 A=C_U(z)=B,\qquad T=C_U(s),\qquad E=T^2.
\]
Projection to $U/B\cong\Z$ is injective on $T$: an element of $T\cap B$
would be fixed by the nontrivial scalar action of $s$.
Thus $T=\langle r\rangle$ is infinite cyclic. Orient $r$ so that its
height $h$ is positive and put $Q=k^{2h}>1$. Definitions
\eqref{arith:eq:I}--\eqref{arith:eq:Div} give, in additive notation,
\[
 I_{r^{2j}}=(Q^{|j|}-1)B\qquad(j\in\Z).
\]
For negative $j$, multiplication by a power of $Q$ is an automorphism
of $B$, so the same subgroup occurs.

Every positive integer $Q^j-1$ is coprime to $k$. For two such integers
$d,e$, one has $eB\subseteq dB$ if and only if $d\mid e$ in $\Z$:
if $e/d\in\Z[1/|k|]$, the primes in $d$ cannot be cancelled by
denominators which are powers of $k$. Lemma~\ref{arith:lem:gcd}, with
$p=1$ and $q=Q$, therefore proves
\eqref{arith:eq:integer-divisibility} in this case as well.

\subsection{The admissible parameters and the interpretation}

To be metabelian is expressed by the sentence
\[
 \operatorname{Met}:=
 \forall x_1,x_2,x_3,x_4\
 [[x_1,x_2],[x_3,x_4]]=1.
\]
For groups in $\mathcal C$, solubility is equivalent to metabelianity,
by Theorem~\ref{arith:thm:soluble}. Define
\[
 \mathcal P_0(z,s):=
 (\neg\operatorname{Met}\land\mathcal P_{\rm n}(z,s))
 \ \lor\
 (\operatorname{Met}\land\mathcal P_{\rm s}(z,s)).
\]
Use the definition of $A$ in \eqref{arith:eq:nonsoluble-A} when
$\neg\operatorname{Met}$ holds, and $A=C_U(z)$ when
$\operatorname{Met}$ holds. This gives one formula $\alpha(x;z,s)$.
The formula defining $E$ is the same in both cases:
\[
 \eta(y;s):=\exists v\,([v,s]=1\land y=v^2).
\]
In \eqref{arith:eq:Div} require also $\eta(y_1;s)$ and $\eta(y_2;s)$.
Finally, put
\begin{align}
 \mathcal P(z,s,u):={}&\mathcal P_0(z,s)\land\eta(u;s)\land u\ne1
 \notag\\
 &{}\land\forall y\,(\eta(y;s)\Rightarrow\Div(u,y;z,s)).
 \label{arith:eq:parameters}
\end{align}

\begin{proof}[Proof of Theorem~\ref{arith:thm:criterion}]
For every pair satisfying $\mathcal P_0$, the calculations above give
$E=\langle r^2\rangle$ and
$\Div(r^{2i},r^{2j};z,s)$ if and only if $i\mid j$.
Consequently, the last two conditions in
\eqref{arith:eq:parameters} select exactly $u=r^2$ and $u=r^{-2}$.
For either choice the map
\[
 \Z\longrightarrow E,\qquad n\longmapsto u^n
\]
identifies $(\Z,+,\mid,0,1)$ with
$(E,\cdot,\Div,1_U,u)$.

The parameter set is nonempty. In the nonsoluble case choose
$z=a^r$ as in Lemma~\ref{arith:lem:find-center} and take an element $s$
with $|\Delta(s)|\ne1$. These satisfy $\mathcal P_{\rm n}$.
In the soluble case choose a nontrivial base element $z$ and a stable
letter $s$; these satisfy $\mathcal P_{\rm s}$. A generator or its
inverse in $E$ supplies $u$.

All the formulas are fixed, equality is ordinary equality, and every
admissible tuple gives the same structure. This proves uniformity,
regularity and injectivity. Applying the classical definition of
multiplication in $(\Z,+,\mid,0,1)$ proves the assertion for the ring.
\end{proof}

\begin{corollary}\label{arith:cor:GBS-unions}
Arithmetic is uniformly, regularly and injectively interpretable in
each of the following classes:
\begin{enumerate}
\item non-unimodular GBS groups;
\item increasing unions of non-unimodular GBS groups whose integral
cohomological dimension is two.
\end{enumerate}
\end{corollary}
\begin{proof}
A GBS group has an aspherical graph of circles of dimension two.
In the non-unimodular case a vertex generator $a$ and an element $s$
satisfy an unequal-power relation, and $\langle a\rangle$ is
commensurated. Apply Theorem~\ref{arith:thm:criterion}.

For $U=\bigcup_iG_i$, choose such $a,s$ in $G_0$. The relation persists
in every $G_i$ and makes $a$ elliptic in its GBS tree. Thus every $G_i$
commensurates $\langle a\rangle$, and so does $U$.
The same theorem applies to the countable group $U$.
\end{proof}

\begin{corollary}\label{arith:thm:main}
The interpretation is uniform over all groups $\BS(m,n)$ with
$m,n\ne0$ and $|m|\ne|n|$.
\end{corollary}

\begin{example}
In $G=\BS(2,3)=\langle a,t\mid ta^2t^{-1}=a^3\rangle$, take
$D=\langle a\rangle=C_G(a)$. The product of the domain and image of
conjugation by $t$ inside $D$ is $D$ itself. On
$E=C_G(t)^2=\langle t^2\rangle$ the construction gives
\[
 I_{t^{2j}}=\langle a^{9^{|j|}-4^{|j|}}\rangle.
\]
For example, $I_{t^2}=\langle a^5\rangle$ and
$I_{t^4}=\langle a^{65}\rangle$.
\end{example}

\subsection{Bi-interpretability}

\begin{definition}
Two structures are \emph{bi-interpretable} if they are interpretable
in one another and both composite interpretations are definably
isomorphic to the original structures. Parameters are allowed in the
interpretations and in the two isomorphisms.
\end{definition}

\begin{proposition}\label{arith:prop:soluble-bi}
For every integer $k$ with $|k|>1$, the group $\BS(1,k)$ is
bi-interpretable with $\Z$, with parameters. Consequently, every
soluble group satisfying the hypotheses of
Theorem~\ref{arith:thm:criterion} is bi-interpretable with $\Z$.
\end{proposition}
\begin{proof}
For $k>1$ the stronger regular bi-interpretability result is
\cite[Theorem 3]{DMBS}. For the remaining case write
\[
 G=\BS(1,-k)=B\rtimes\langle t\rangle,\qquad
 B=\Z[1/k],\quad k>1,
\]
where $t$ acts on $B$ by multiplication by $-k$. Let $a$ be the
standard base generator, put $s=t^2$, and set
\[
 H=B\rtimes\langle s\rangle\cong\BS(1,k^2).
\]
Since $B=C_G(a)$ and $C_G(t)=\langle t\rangle$, the subgroup
$H=C_G(a)C_G(t)^2$ is definable in $G$.

Inside $H$, one has $B=C_H(a)$ and $\langle s\rangle=C_H(s)$.
Every $h\in H$ has a unique expression $h=bv$ with $b\in B$ and
$v\in\langle s\rangle$. The restriction of conjugation by $t$ is
therefore the definable automorphism
\[
 \tau(bv)=b^{-k}v.
\]
On $H\times\{0,1\}$ define
\[
 (h,i)(h',j)=
 \bigl(h\,\tau^i(h')s^{ij},\,i+j\pmod2\bigr).
\]
The two-element set is coded by $1_H,s$, and this multiplication is
a finite list of first-order formulas. The map $(h,i)\mapsto ht^i$
is an isomorphism onto $G$ and is definable in $G$.
Applying the definition of $H$ in this interpreted group gives
$H\times\{0\}$, with the definable isomorphism $h\mapsto(h,0)$.
Thus $G$ and $H$ are bi-interpretable. Compose with the known
bi-interpretation between $H\cong\BS(1,k^2)$ and $\Z$.

The final assertion follows from Theorem~\ref{arith:thm:soluble} and
the exclusion of $k=1,-1$ by the unequal-power relation.
\end{proof}

For fixed $k>1$, the construction for $\BS(1,-k)$ can also be made
regular. The condition
\[
 a\ne1,\qquad tat^{-1}=a^{-k}
\]
forces $a$ to be a nontrivial base element and $t$ to have height one.
Hence the displayed definitions of $H$ and $\tau$ work for every such
pair. Append the admissible parameters from the regular
bi-interpretation of $H$, with their defining condition relativized
to $H$. The two displayed isomorphisms are defined by the same
formulas for every admissible tuple.

\begin{proposition}\label{arith:prop:no-general-bi}
There is a group in $\mathcal C$ which is not bi-interpretable with
$\Z$, even with finitely many parameters.
\end{proposition}
\begin{proof}
Let
\[
 U=\langle a,t_0,t_1,\ldots\mid
                 t_iat_i^{-1}=a^2\ (i\geq0)\rangle.
\]
This is the fundamental group of a graph of infinite cyclic groups
with one vertex and countably many loop edges. Its graph of circles
is an aspherical two-dimensional complex. It contains $\BS(1,2)$,
so $\cdim U=2$. The subgroup $\langle a\rangle$ is commensurated and
$\Delta(t_i)=2$ for all $i$.

Given a finite tuple in $U$, choose words representing its entries.
Only finitely many $t_i$ occur. Interchanging two unused stable
letters and fixing the remaining generators gives a nontrivial
automorphism fixing the tuple.

If $U$ were bi-interpretable with $\Z$, include in one finite tuple
all parameters used in the interpretations and their composite
isomorphisms. An automorphism fixing this tuple induces an automorphism
of the interpreted integer ring, which must be the identity.
The definable isomorphism for the composite interpretation then
forces the original automorphism of $U$ to be the identity.
This contradicts the preceding construction.
\end{proof}


\section{Commensurators and arithmetic in one-relator groups}
\label{sec:nonuni-ambient}

\begin{theorem}\label{ambient:thm:main}
The ring $(\Z,+,\times,0,1)$ is uniformly, regularly and injectively
interpretable in the class of one-relator groups containing a
non-unimodular Baumslag--Solitar subgroup. The interpretation formulas
and the formula defining admissible parameters are independent of
the Baumslag--Solitar exponents.
\end{theorem}

A one-relator group $G$ containing such a subgroup is torsion-free:
one-relator groups with torsion are hyperbolic, whereas unequal powers
of an infinite-order element of a hyperbolic group cannot be conjugate.
Lyndon's asphericity theorem gives $\cdim G\leq2$ \cite{Lyndon}, and
this bound passes to subgroups \cite[Chapter VIII]{Brown}.

\subsection{The commensurator of the base cyclic subgroup}

Fix
\[
 H=\langle a,s\rangle\cong\BS(m,n),\qquad
 sa^ms^{-1}=a^n,\quad |m|\ne|n|,
\]
and put $M=\Comm_G(\langle a\rangle)$. The group $M$ is countable,
has cohomological dimension two, and commensurates $\langle a\rangle$.
It therefore satisfies the hypotheses of
Theorem~\ref{arith:thm:criterion}.

\begin{lemma}\label{ambient:lem:intrinsic}
For every nonzero integer $r$,
\[
 C_G(a^r)\leq M,\qquad C_G(a^r)=C_M(a^r).
\]
\end{lemma}
\begin{proof}
An element centralizing $a^r$ conjugates $\langle a\rangle$ to an
infinite cyclic subgroup containing $a^r$. Their intersection is
nontrivial, so the element belongs to $M$.
\end{proof}

\begin{proposition}\label{ambient:prop:defM}
The subgroup $M$ is definable in $G$ with parameters. More precisely,
one of the following holds.
\begin{enumerate}
\item Some $z=a^r$, $r\geq1$, has nonabelian centralizer. Then
$D=Z(C_G(z))$ is infinite cyclic and
\begin{equation}\label{ambient:eq:mu-n}
 M=\{g\in G:\exists x\,
      (\delta(x;z)\land x\ne1\land\delta(gxg^{-1};z))\}.
\end{equation}
\item Every $C_G(a^r)$, $r\geq1$, is abelian. Then
$M\cong\BS(1,k)$ for some $|k|>1$, and
\begin{equation}\label{ambient:eq:mu-s}
 M=N_G(C_G(a)).
\end{equation}
\end{enumerate}
\end{proposition}
\begin{proof}
In the first case, Lemma~\ref{ambient:lem:center} makes $D$ cyclic.
It contains $a^r$, so it is commensurable with $\langle a\rangle$.
Consequently $\Comm_G(D)=M$, and nontrivial intersection of cyclic
groups gives \eqref{ambient:eq:mu-n}.

In the second case, Lemmas~\ref{arith:lem:modular} and
\ref{ambient:lem:intrinsic} show that
$\ker(\Delta:M\to\Q^\times)$ is abelian. Thus $M$ is metabelian.
Theorem~\ref{arith:thm:soluble} and the unequal-power relation give
$M\cong\BS(1,k)$, $|k|>1$.

Let $B\cong\Z[1/|k|]$ be its standard base. The image of $a$ in
$M/B\cong\Z$ is zero, since the displayed relation gives
$m\bar a=n\bar a$ and $m\ne n$. Thus $C_M(a)=B$, and
Lemma~\ref{ambient:lem:intrinsic} gives $C_G(a)=B$.
Every element of $M$ normalizes $B$. Conversely, if $g\in N_G(B)$,
then $a$ and $gag^{-1}$ are nontrivial elements of the torsion-free
locally cyclic group $B$. They have equal nonzero powers, so $g\in M$.
This proves \eqref{ambient:eq:mu-s}.
\end{proof}

\begin{proposition}\label{ambient:prop:GBS-overgroups}
Every finitely generated subgroup $J$ with $H\leq J\leq M$ is
non-unimodular GBS. Every GBS subgroup of $G$ containing $H$ is
contained in $M$. Consequently $M$ is the union of an increasing
sequence of non-unimodular GBS overgroups of $H$.

Moreover, $M$ is finitely generated if and only if there is an
inclusion-maximal GBS subgroup of $G$ containing $H$. In this event
the maximal subgroup is $M$, and it contains every GBS overgroup
of $H$ in $G$.
\end{proposition}
\begin{proof}
The first assertion is Lemma~\ref{arith:lem:local} applied inside $M$.
If $K$ is GBS and contains $H$, the unequal-power relation makes $a$
elliptic in its Bass--Serre tree. All of $K$ therefore commensurates
$\langle a\rangle$, and $K\leq M$. Enumerate $M=\{g_1,g_2,\ldots\}$;
the groups $\langle H,g_1,\ldots,g_i\rangle$ give the required chain.

If $M$ is finitely generated, the first assertion makes it GBS, and
the second makes it the largest such overgroup. Conversely, let $K$
be an inclusion-maximal GBS overgroup of $H$. For every $x\in M$,
the group $\langle K,x\rangle$ is a finitely generated intermediate
subgroup and hence GBS. Maximality gives $x\in K$. Thus $K=M$.
\end{proof}

All applications of Kropholler's theorem here are to finitely
generated subgroups. We do not assume that the commensurator is
finitely generated.

\subsection{A definable family of admissible subgroups}

To prove regularity in Theorem~\ref{ambient:thm:main}, we give
first-order conditions on the parameters defining the intermediate
subgroup. Write
\[
 \operatorname{AbC}(z):=
 \forall x\,\forall y\,
 (([x,z]=1\land[y,z]=1)\Rightarrow[x,y]=1).
\]
For $z\in G$, define
\begin{align*}
 \mu_{\rm n}(g;z)&:=
 \exists x\,(\delta(x;z)\land x\ne1\land\delta(gxg^{-1};z)),\\
 \mu_{\rm s}(g;z)&:=
 \forall x\,([x,z]=1\ \Longleftrightarrow\ [gxg^{-1},z]=1).
\end{align*}
The second formula defines $N_z=N_G(C_G(z))$.
When $z\ne1$ and $\neg\operatorname{AbC}(z)$, the first defines
$\Comm_G(D_z)$, where $D_z=Z(C_G(z))$ is infinite cyclic.
Put
\[
 R_s(z)=\{x\in D_z:sxs^{-1}\in D_z\},\qquad
 S_s(z)=\{x\in D_z:s^{-1}xs\in D_z\}.
\]

Define the first parameter condition by
\[
 \mathcal Q_{\rm n}(z):=
 z\ne1\land\neg\operatorname{AbC}(z)\land
 \exists s\,\bigl(\mu_{\rm n}(s;z)\land R_s(z)\ne S_s(z)\bigr).
\]
For the second condition, relativize the metabelian identity to $N_z$:
\[
 \mathcal Q_{\rm s}(z):=
 z\ne1\land\operatorname{AbC}(z)\land
 \operatorname{Met}^{N_z}\land
 \exists x,y\in N_z\ ([x^2,y^2]\ne1).
\]
All subgroup-membership conditions and restricted quantifiers in
these expressions abbreviate the displayed first-order formulas.
Let
\begin{align}
 \mathcal Q(z)&:=\mathcal Q_{\rm n}(z)\lor\mathcal Q_{\rm s}(z),
 \label{ambient:eq:Q}\\
 \mu(g;z)&:=
 (\neg\operatorname{AbC}(z)\land\mu_{\rm n}(g;z))
 \lor(\operatorname{AbC}(z)\land\mu_{\rm s}(g;z)).
 \label{ambient:eq:mu}
\end{align}

\begin{lemma}\label{ambient:lem:admissible-M}
If $\mathcal Q(z)$ holds, then $M_z=\mu(G;z)$ belongs to $\mathcal C$.
If $G$ contains a non-unimodular Baumslag--Solitar subgroup, then
$\mathcal Q(G)$ is nonempty.
\end{lemma}
\begin{proof}
Suppose $\mathcal Q_{\rm n}(z)$ holds. The group $M_z$ commensurates
the infinite cyclic group $D_z$. The inequality $R_s(z)\ne S_s(z)$
says that conjugation by $s$ has modular value of absolute value
different from one. Thus $M_z$ satisfies the hypotheses of
Theorem~\ref{arith:thm:criterion}.

Suppose $\mathcal Q_{\rm s}(z)$ holds. The group $M_z=N_z$ is a
countable metabelian subgroup of $G$ and has cohomological dimension
at most two. By Theorem~\ref{arith:thm:soluble}, it is locally cyclic
or isomorphic to $\BS(1,k)$. Locally cyclic groups and $\BS(1,1)$
are abelian. In $\BS(1,-1)$, all squares lie in the normal abelian
subgroup of index two, so they commute. The last condition in
$\mathcal Q_{\rm s}$ excludes these cases. Hence
$M_z\cong\BS(1,k)$ with $|k|>1$, and $M_z\in\mathcal C$.

For nonemptiness, apply Proposition~\ref{ambient:prop:defM} to the
given subgroup $H$. In its first case choose the indicated power
$z=a^r$. The modular value of the original stable letter supplies
$\mathcal Q_{\rm n}(z)$. In its second case take $z=a$.
Then $N_z=M\cong\BS(1,k)$, $|k|>1$, and the squares of a nontrivial
base element and a stable letter do not commute. Thus
$\mathcal Q_{\rm s}(z)$ holds.
\end{proof}

\begin{proof}[Proof of Theorem~\ref{ambient:thm:main}]
For every $w$ satisfying $\mathcal Q(w)$, apply
Theorem~\ref{arith:thm:criterion} inside $M_w$.
Let $\mathcal P$ be the parameter formula
\eqref{arith:eq:parameters}. A parameter formula in $G$ is
\[
 \mathcal P_G(w,z,s,u):=
 \mathcal Q(w)\land
 \mu(z;w)\land\mu(s;w)\land\mu(u;w)
 \land\mathcal P^{M_w}(z,s,u),
\]
where every quantifier in $\mathcal P^{M_w}$ is restricted to
$\mu(G;w)$. Relativize the interpretation formulas in the same way.

Lemma~\ref{ambient:lem:admissible-M} and
Theorem~\ref{arith:thm:criterion} show that $\mathcal P_G(G)$ is
nonempty and that every tuple satisfying it defines the integer ring.
Equality remains ordinary equality. None of these formulas contains
the exponents of the originally given Baumslag--Solitar subgroup.
This proves all three assertions.
\end{proof}

\begin{corollary}\label{ambient:cor:intrinsic-arithmetic}
With the notation above, arithmetic is uniformly, regularly and
injectively interpretable in the pure group $M$.
If $M$ is soluble, it is bi-interpretable with $\Z$, with parameters.
\end{corollary}
\begin{proof}
Apply Theorem~\ref{arith:thm:criterion} and
Proposition~\ref{arith:prop:soluble-bi}, respectively, to $M$.
\end{proof}

The second assertion is a statement about the commensurator $M$.
Bi-interpretability of the ambient group $G$ requires definable
isomorphisms for the two composite interpretations.

\subsection{Further interpretations and undecidability}

\begin{lemma}\label{ambient:lem:coding}
Every fixed finitely presented group is absolutely and injectively
interpretable in $(\Z,+,\times,0,1)$.
\end{lemma}
\begin{proof}
Fix $H=\langle x_1,\ldots,x_r\mid R_1,\ldots,R_k\rangle$ and code
finite words in $x_1^{\pm1},\ldots,x_r^{\pm1}$ by natural numbers.
The natural numbers are definable in $\Z$ as sums of four squares.
The set $W$ of word codes, concatenation and inversion are computable,
hence arithmetically definable.

The relation $F(u,v)$ saying that two words represent the same
element of $H$ is computably enumerable: enumerate applications
of the defining relators and free reductions. It is therefore
arithmetically definable; see \cite[Section 2.1]{DMBS}.
For an injective interpretation, use
\[
 W_{\min}=\{u\in W:\forall v\in W\,
                    (v<u\Rightarrow\neg F(v,u))\}.
\]
Each equivalence class has one least code. On $W_{\min}$ define
multiplication by $F(\operatorname{concat}(u,v),w)$, and define
inversion and identity in the same way. Equality is ordinary
equality. The formulas depend only on the fixed finite presentation.
\end{proof}

\begin{corollary}\label{ambient:cor:fp}
Every fixed finitely presented group is uniformly, regularly and
injectively interpretable in every group in $\mathcal C$, and in
every one-relator group containing a non-unimodular Baumslag--Solitar
subgroup.
\end{corollary}
\begin{proof}
Compose Lemma~\ref{ambient:lem:coding} with
Theorem~\ref{arith:thm:criterion} or
Theorem~\ref{ambient:thm:main}. The admissible parameter set is unchanged.
\end{proof}

These group operations are obtained by arithmetic coding. The
definable subgroup in Theorem~\ref{thm:bs-interpretation}, in contrast,
has the inherited group operations.

\begin{corollary}\label{ambient:cor:undecidable}
\label{cor:all-nonuni-containing}
Every one-relator group containing a non-unimodular Baumslag--Solitar
subgroup has undecidable elementary theory. One computable translation
of arithmetic sentences works for all groups in this class.
The same assertion holds for the class $\mathcal C$.
\end{corollary}
\begin{proof}
Let $P(\bar c)$ be the relevant admissibility formula and let
$\sigma^I(\bar c)$ be the translation of an arithmetic sentence
$\sigma$ under the interpretation. Put
\[
 \widehat{\sigma}:=
 \exists\bar c\,(P(\bar c)\land\sigma^I(\bar c)).
\]
Every admissible tuple gives the standard integer ring, and there
is at least one such tuple. Hence, for every group $V$ in the stated
class,
\[
 V\models\widehat{\sigma}
 \quad\Longleftrightarrow\quad
 \Z\models\sigma.
\]
The translation is computable and independent of $V$.
Undecidability follows from that of arithmetic \cite{TMR}.
\end{proof}

The results of this section do not use either result from
\cite{OAI,OAISpecial}. Their defining formulas
contain universal quantifiers; the Diophantine problems are considered
separately in Section~\ref{sec:beyond}.

\section{Equations and other decision problems}\label{sec:beyond}

\begin{proposition}\label{prop:bs-theory}
For every integer \(n>1\), \(\BS(1,n)\) has undecidable elementary
theory, but is equationally Noetherian.
\end{proposition}
\begin{proof}
The arithmetic interpretation is due to Khelif \cite{Khelif};
the stronger regular bi-interpretation and elementary-equivalence
classification are given in \cite{KMS}, \cite{DMBS}.
For Noetherianity, the standard faithful representation is
\[
a\longmapsto\begin{pmatrix}1&1\\0&1\end{pmatrix},\qquad
t\longmapsto\begin{pmatrix}n&0\\0&1\end{pmatrix}.
\]
It identifies \(\Z[1/n]\rtimes\Z\) with a subgroup of
\(\mathrm{GL}_2(\Q)\). Linear groups over a field are equationally
Noetherian: entries of equations, with inverse determinants added
as variables, generate an ideal in a Noetherian polynomial ring.
Finite generation of that ideal gives a finite subsystem.
\end{proof}

\begin{proposition}\label{prop:en}
For nonzero integers \(m,n\), \(\BS(m,n)\) is equationally Noetherian
if and only if
\[
\min\{|m|,|n|\}=1\quad\text{or}\quad |m|=|n|.
\]
Moreover, \(\BS(m,\pm m)\) has decidable equations and inequations.
\end{proposition}
\begin{proof}
The Noetherianity criterion is
\cite[Proposition 5.18]{GuyotStalder}, citing \cite{BMR}.
The unimodular groups have the explicit finite-index direct products
of Proposition~\ref{prop:unimodular-products}; apply
\cite[Theorem 3.1]{CHR} for equations and inequations, including
recognizable constraints.
\end{proof}

In particular, the groups \(\BS(m,\pm m)\) are equationally
Noetherian. 

\begin{theorem}\label{thm:bg}
The Diophantine problem is undecidable in the one-relator group
\[
B_d=\langle a,t\mid(tat^{-1})a(tat^{-1})^{-1}=a^d\rangle
\]
for every integer \(d\geq2\).
\end{theorem}

For prime \(d\), this is Gray--Levine
\cite[Theorem 1.1]{GrayLevine}; Nyberg-Brodda
\cite[Theorem A]{NybergBG} gives all \(d\geq2\).
Thus the Diophantine problem is not decidable for all one-relator
groups. Each \(B_d\) contains \(\BS(1,d)\), by
Lemma~\ref{lem:power}, taking \(u=a\) and \(v=tat^{-1}\).
Corollary~\ref{cor:containing-undecidability} also gives undecidability
of \(\Th(B_d)\).

This does not settle arbitrary systems of equations in \(\BS(1,n)\)
itself. 
Decidability results for quadratic equations, for example
\cite{MU}, concern a restricted class of equations. Likewise,
undecidability of a full elementary theory does not imply
undecidability of its positive existential fragment.

\begin{proposition}\label{prop:limit}
There is an unconditional algorithm deciding whether an arbitrary
one-relator group is fully residually free.
\end{proposition}
\begin{proof}
A finitely generated group is \emph{fully residually free} if every
finite set of nonidentity elements survives simultaneously under
some homomorphism to a free group; these are the limit groups. Bumagin- Kharlampovich - Miasnikov \cite{BKM} and 
Groves--Wilton \cite{GWLimit} gives a
recognition algorithm from a finite presentation and a solution to
its word problem. Magnus supplies the latter uniformly for one-relator
groups. 
\end{proof}

\begin{theorem}\label{thm:commutator}
For \(m\geq1\), the one-relator group
\[
K_m=\langle a,t\mid[t,a(a^{-1}t^{-1}at)^{-m}]=1\rangle
\]
has a retract isomorphic to \(\BS(m,m+1)\).
The group \(K_2\) is non-Hopfian and not residually finite,
and \(K_1\) is not automatic.
\end{theorem}
This is \cite[Theorem 1.3]{WangZhang}. In particular, a one-relator
group with a commutator relator need not be Hopfian, residually finite,
or automatic. Since \(K_1\) contains \(\BS(1,2)\),
Corollary~\ref{cor:containing-undecidability} applies to \(K_1\).
For every \(m\ge2\), Corollary~\ref{arith:thm:main} and
Corollary~\ref{ambient:cor:undecidable} also imply that
\(\Th(K_m)\) is undecidable, since the retract embeds in \(K_m\).
Gardam--Woodhouse \cite{GardamWoodhouse} give further nonautomatic
one-relator groups containing no \(\BS(m,n)\) with \(m\neq\pm n\).

\section{Unimodular Baumslag--Solitar groups}\label{sec:unimodular}

By changing the sign of both parameters, every Baumslag--Solitar group
with \(|m|=|n|=k>1\) is isomorphic to \(\BS(k,k)\) or \(\BS(k,-k)\).
Write \(F_r\) for the free group of rank \(r\).

\begin{proposition}\label{prop:unimodular-products}
For \(k>1\), the following normal finite-index subgroups are direct products:
\begin{enumerate}
\item In \(G_+=\BS(k,k)\), the kernel of
\[
 G_+\longrightarrow C_k,\qquad a\longmapsto\overline1,\quad
 t\longmapsto\overline0
\]
has index \(k\) and is isomorphic to \(F_k\times\Z\).
\item In \(G_-=\BS(k,-k)\), the kernel of
\[
 G_-\longrightarrow C_k\times C_2,\qquad
 a\longmapsto(\overline1,\overline0),\quad
 t\longmapsto(\overline0,\overline1)
\]
has index \(2k\) and is isomorphic to \(F_{2k-1}\times\Z\).
\end{enumerate}
\end{proposition}
\begin{proof}
For the first assertion, put
\[
 z=a^k,\qquad x_i=a^it a^{-i}\quad(0\le i<k).
\]
Reidemeister--Schreier rewriting, with representatives
\(1,a,\ldots,a^{k-1}\), gives the presentation
\[
 H_+=\langle z,x_0,\ldots,x_{k-1}\mid[z,x_i]=1\ (0\le i<k)\rangle.
\]
Thus \(H_+=F(x_0,\ldots,x_{k-1})\times\langle z\rangle\), and the map
onto \(C_k\) gives its index.

For the second assertion, let \(A=\langle a^k\rangle\). This is an
infinite cyclic normal subgroup of \(G_-\); conjugation by \(a\) fixes its
generator and conjugation by \(t\) inverts it. The quotient is
\[
 G_-/A\cong C_k*\Z.
\]
Let \(H_-\) be the kernel in the statement and let \(\overline H_-\) be
its image in this quotient. Then \(\overline H_-\) is the kernel of the
corresponding surjection \(C_k*\Z\to C_k\times C_2\).

View \(C_k*\Z\) as the fundamental group of a graph of groups with one
vertex group \(C_k\) and one loop edge with trivial group. The restriction
of the finite quotient map to every conjugate of \(C_k\) is injective, so
\(\overline H_-\) acts freely on the Bass--Serre tree. Its quotient graph
has \(2k/k=2\) vertices and \(2k\) edges. It follows that
\(\overline H_-\cong F_{2k-1}\).

Every element of \(H_-\) has even exponent sum in \(t\), hence centralizes
\(A\). The sequence
\[
 1\longrightarrow A\longrightarrow H_-\longrightarrow F_{2k-1}
 \longrightarrow1
\]
is therefore central and splits: choose lifts of a free basis of the
quotient. Their generated subgroup maps isomorphically onto the free
quotient and meets \(A\) trivially. Hence \(H_-=F_{2k-1}\times A\).
\end{proof}

\begin{example}
The group \(\BS(2,2)\) has the index-two subgroup
\[
 \langle t,ata^{-1},a^2\rangle
 =F(t,ata^{-1})\times\langle a^2\rangle.
\]
The group \(\BS(2,-2)\) has an index-four subgroup isomorphic to
\(F_3\times\Z\).
\end{example}

\subsection{The elementary theories of the direct products}

We record the effective separation of formulas which will also be used
for the periodic action. It is the finite-product instance of the
Feferman--Vaught method \cite{FV}.

\begin{lemma}\label{lem:sort-separation}
Let \(\mathcal A,\mathcal B\) be structures in disjoint languages,
regarded as a two-sorted structure with no symbols connecting the sorts.
Every formula is effectively equivalent to a finite disjunction
\[
 \bigvee_{\nu=1}^{N}
 \bigl(\alpha_\nu(\bar x)\land\beta_\nu(\bar y)\bigr),
\]
where \(\alpha_\nu\) uses only the \(\mathcal A\)-sort and
\(\beta_\nu\) uses only the \(\mathcal B\)-sort. If the theories of both
sorts are decidable, then the theory of the two-sorted structure is
decidable.
\end{lemma}
\begin{proof}
Atomic formulas belong to one sort. Finite disjunctions of separated
conjunctions are closed effectively under Boolean operations, by
distributivity and disjunctive normal form. For an \(\mathcal A\)-variable
\(x\), existential quantification gives
\[
 \exists x\,\bigvee_\nu(\alpha_\nu\land\beta_\nu)
 \ \longleftrightarrow\
 \bigvee_\nu((\exists x\,\alpha_\nu)\land\beta_\nu),
\]
because \(x\) does not occur in \(\beta_\nu\). The analogous operation
handles the other sort. Universal quantification follows from negation
and existential quantification. Induction on formulas proves the
assertion, including arbitrary quantifier alternation. For a sentence,
the resulting \(\alpha_\nu,\beta_\nu\) are sentences in their respective
sorts and their truth values can be decided.
\end{proof}

\begin{corollary}\label{cor:free-times-z}
For every fixed \(r\ge1\), the elementary theory of \(F_r\times\Z\)
is decidable. In particular, both direct-product subgroups in
Proposition~\ref{prop:unimodular-products} have decidable elementary theory.
\end{corollary}
\begin{proof}
The theories of \(F_r\) and of \((\Z,+)\) are decidable, by
\cite{KMFree} and Presburger arithmetic; the rank-one case is also
additive arithmetic. Interpret the product on pairs with coordinatewise
multiplication and inverse. Each atomic equality becomes a conjunction
of one equality in each sort. Apply Lemma~\ref{lem:sort-separation}.
\end{proof}

\subsection{The periodic action and the full theory}

\begin{proposition}\label{prop:periodic-action}
Let \(F_k=F(x_0,\ldots,x_{k-1})\) and let
\(\sigma\in\operatorname{Aut}(F_k)\) be defined by
\(\sigma(x_i)=x_{i+1\pmod k}\). Then
\[
 \BS(k,k)\cong F_k\rtimes_\sigma\Z,
\]
where the generator of \(\Z\) corresponds to \(a\), and
\(x_i\) corresponds to \(a^it a^{-i}\).
\end{proposition}
\begin{proof}
The semidirect product has presentation
\[
 \langle a,x_0,\ldots,x_{k-1}\mid
 ax_i a^{-1}=x_{i+1\pmod k}\ (0\le i<k)\rangle.
\]
Eliminating \(x_1,\ldots,x_{k-1}\) gives
\(\langle a,x_0\mid a^kx_0a^{-k}=x_0\rangle\).
Taking \(x_0=t\) yields the defining presentation of \(\BS(k,k)\).
\end{proof}

\begin{proposition}\label{prop:periodic-reduction}
If \(\Th(F_k,\cdot,{}^{-1},1,\sigma)\) is decidable for the automorphism
of Proposition~\ref{prop:periodic-action}, then
\(\Th(\BS(k,k))\) is decidable. There is an effective translation of
every group sentence into a finite Boolean combination of sentences in
\((F_k,\sigma)\) and \((\Z,+,0,1)\).
\end{proposition}
\begin{proof}
Every element has a unique expression \(wa^{kh+r}\), where
\(w\in F_k\), \(h\in\Z\), and \(0\le r<k\). Use the coordinates
\((w,h,r)\). Since \(\sigma^k=1\), multiplication is
\begin{equation}\label{eq:periodic-mult}
 (w,h,r)(v,j,s)=
 \left(w\sigma^r(v),\ h+j+\left\lfloor\frac{r+s}{k}\right\rfloor,
 (r+s)\bmod k\right).
\end{equation}
Inversion is
\begin{equation}\label{eq:periodic-inverse}
 (w,h,r)^{-1}=
 \begin{cases}
 (w^{-1},-h,0),&r=0,\\
 (\sigma^{-r}(w^{-1}),-h-1,k-r),&r>0.
 \end{cases}
\end{equation}
For fixed residue coordinates, every group term has a free-group
coordinate which is a term in the language with \(\sigma\), and an
integer coordinate which is an integral linear expression. An atomic
equality is false when its residue coordinates differ; otherwise it is
a conjunction of an equality in the free-group sort and a linear
equality in the integer sort.

Translate quantifiers by the finite residue cases:
\[
 \exists g\,\varphi(g)
 \ \longleftrightarrow\
 \bigvee_{r=0}^{k-1}\exists w\,\exists h\,
 \widehat\varphi(w,h,r).
\]
Universal quantifiers are treated by conjunction over the residues and
universal quantification over the two infinite coordinates. Apply
Lemma~\ref{lem:sort-separation}. Presburger arithmetic decides the
integer sentences, and the assumed algorithm decides the free-group
sentences. Thus the procedure terminates on every input sentence.
\end{proof}

\begin{remark}
Proposition~\ref{prop:periodic-reduction} is conditional. Decidability of
the pure free-group theory does not supply a decision procedure for its
expansion by \(\sigma\). The images of the basis determine a computable
automorphism, but a first-order reduction also has to retain the relation
between a quantified element and its image. Replacing
\(w,\sigma(w),\ldots,\sigma^{k-1}(w)\) by independent variables loses
that relation. For example, \(\sigma\) fixes no nonidentity element: a
nonempty reduced word changes its first letter under the cyclic basis
permutation. Hence
\[
 \exists w\,(w\ne1\land\sigma(w)=w)
\]
is false, whereas independent orbit coordinates would permit a
nonidentity solution to the corresponding equality.
\end{remark}

There is a second formulation of the sufficient condition. Put
\[
 V_k=F_k\rtimes_\sigma C_k\cong C_k*\Z,
\]
where the generator \(c\) of \(C_k\) conjugates \(x_i\) to
\(x_{i+1\pmod k}\). Eliminating the \(x_i\), \(i>0\), gives the
presentation \(\langle c,x_0\mid c^k=1\rangle\). Let
\[
 \rho:V_k\longrightarrow C_k,\qquad
 \rho(c)=c,\quad\rho(x_0)=1.
\]
Then \(\ker\rho=F_k\), and \(\sigma\) is conjugation by \(c\) on
this kernel. Consequently, decidability of
\(\Th(V_k,\ker\rho,c)\) would imply decidability of \(\Th(F_k,\sigma)\)
and hence of \(\Th(\BS(k,k))\). Here \(\ker\rho\) is a named unary
predicate; its definability in the pure group is not assumed.

The analysis of free groups expanded by finite-quotient predicates in
\cite{deLaNuez} supplies related Makanin--Razborov diagrams, towers and
formal-solution results. Those results on positive theory do not give
the full decision theorem required in Proposition~\ref{prop:periodic-reduction}.
For \(\BS(k,-k)\), Proposition~\ref{prop:unimodular-products} proves the
virtual direct-product statement. 

\section{Effective virtual free-by-cyclic presentations}\label{sec:effective-vfc}

Assuming the virtual compact specialness result stated as Theorem~1.1
in \cite{OAISpecial}, we construct a virtual free-by-cyclic presentation for
a hyperbolic one-relator group. We compute a finite-index normal subgroup,
an embedding in a mapping torus of a finite-rank free-group automorphism,
and a finite HNN presentation. The free kernel may have infinite rank.

\subsection{Virtual free-by-cyclic subgroups}

\begin{definition}
A group $H$ is \emph{free-by-cyclic} if there is an exact sequence
\[
1\longrightarrow K\longrightarrow H\longrightarrow C\longrightarrow 1,
\]
where $K$ is free and $C\leq\Z$. The rank of $K$ is not required to be
finite, and $C$ may be trivial. A group is \emph{virtually free-by-cyclic}
if it has a finite-index subgroup with this property.
\end{definition}

If $C\neq 0$, identifying $C$ with $\Z$ and choosing a preimage $s$ of
$1$ gives $H=K\rtimes_\beta\langle s\rangle$, where
$\beta(k)=sks^{-1}$. We use the convention
\begin{equation}\label{vfc:eq:ambient}
E_\alpha=F_d\rtimes_\alpha\langle t\rangle
=\left\langle a_1,\ldots,a_d,t\ \middle|\
ta_it^{-1}=\alpha(a_i),\ 1\leq i\leq d\right\rangle.
\end{equation}
Thus elements of $E_\alpha$ have normal forms $wt^n$, with $w\in F_d$.
The projection $\pi:E_\alpha\to\Z$ is given by $\pi(wt^n)=n$.

The result stated as Theorem~1.1 in \cite{OAISpecial} asserts that every
hyperbolic one-relator group is virtually compact special.
Corollary~18.1 of that paper then
gives a finite-index subgroup embedded in a group of the form
\eqref{vfc:eq:ambient}. The group $F_d$ has finite rank, but the intersection of the
embedded subgroup with $F_d$ need not be finitely generated.
As explained in the introduction, we have read that paper but have not
verified all details of its proofs. The following theorem assumes the
cited virtual compact specialness result.

\begin{theorem}\label{vfc:thm:main}
Assume the result stated as Theorem~1.1 in \cite{OAISpecial}.
There is a uniform algorithm which,
given a finite one-relator presentation
\[
G=\langle X\mid r\rangle
\]
of a hyperbolic group, computes the following data:
\begin{enumerate}[label=(\roman*)]
\item a finite-index normal subgroup $H\trianglelefteq G$, a finite
quotient with kernel $H$, generators of $H$ as words in $X$, and a
finite presentation of $H$;
\item an integer $d\geq0$, an automorphism $\alpha\in\Aut(F_d)$,
specified together with its inverse, and an embedding
$\iota:H\hookrightarrow E_\alpha$, specified on the generators of $H$;
\item a homomorphism $\chi:H\to C$, where $C=0$ or $C=\Z$, with free
kernel $K$; if $C=\Z$, the algorithm also gives an element $s$ with
$\chi(s)=1$;
\item if $C=0$, a finite free basis of $H$; if $C=\Z$, a finite HNN
presentation of $H$ with finitely generated free base and associated
subgroups, and an effective free basis of $K$ on which conjugation by
$s$ and by $s^{-1}$ can be computed.
\end{enumerate}
An effective basis in (iv) means a basis indexed by a decidable set,
with a computable word for each basis element. The set may be infinite.
The images under conjugation are computed as finite words in the basis.
\end{theorem}

To prove termination, we shall use the following quasiconvex
embedding theorem.

\begin{lemma}\label{vfc:lem:structural}
Let $G$ be a hyperbolic virtually compact special group which is
virtually free-by-cyclic. Then either $G$ is virtually free, or a
finite-index normal subgroup of $G$ embeds quasiconvexly in a
hyperbolic group $E_\alpha$ of the form \eqref{vfc:eq:ambient}.
\end{lemma}

\begin{proof}
Theorem~1.1 of Kielak--Linton \cite{KL} gives
$\cdq(G)\leq2$ and $b_2^{(2)}(G)=0$. If $\cdq(G)\leq1$, the group is
virtually free. In dimension two, their Theorem~1.12 and the proof of
Theorem~1.1 give a finite-index subgroup $L\leq G$ embedded
quasiconvexly in a hyperbolic group with a finitely generated free
normal subgroup and quotient $\Z$. Splitting the quotient gives
\eqref{vfc:eq:ambient}. The normal core of $L$ in $G$ has finite index in
$L$, so its image remains quasiconvex.
\end{proof}

For a hyperbolic one-relator group, the virtual compact specialness
result stated in \cite{OAISpecial}, together with \cite[Corollary~1.3]{KL}, gives
the hypotheses of Lemma~\ref{vfc:lem:structural}. This also applies to
trivial and proper-power relators. If $G$ is finite, the trivial subgroup
gives the virtually free case.

\subsection{Construction of an embedding}

From a finite presentation of a hyperbolic group one can compute a
hyperbolicity constant, a solution of the word problem, and a shortlex
automatic structure; see \cite[Section~2.4]{KMW}. The procedure detecting
hyperbolicity terminates precisely on hyperbolic presentations.

\begin{lemma}\label{vfc:lem:qc}
There is a uniform partial algorithm which, given a finite presentation
of a hyperbolic group $E$ and words $y_1,\ldots,y_k$ in its generators,
terminates precisely when $Q=\langle y_1,\ldots,y_k\rangle$ is
quasiconvex. When it terminates, it computes a membership algorithm for
$Q$ and a finite presentation
\[
Q=\langle y_1,\ldots,y_k\mid S\rangle
\]
with the indicated identification of generators.
\end{lemma}

\begin{proof}
Kapovich's algorithm \cite{Kap} detects quasiconvexity. Alternatively,
\cite[Theorem~4.1 and Corollary~4.2]{KMW} compute the Stallings graph
and a quasiconvexity constant. These data solve membership in $Q$.

They also give an effective quasi-isometric embedding of the Cayley
graph of $Q$, with the indicated finite generating set, in that of $E$.
The effective thin-triangle and Morse bounds therefore give a
hyperbolicity constant for $Q$. The usual presentation by sufficiently
short trivial words in a hyperbolic group now gives $S$: enumerate the
words up to the resulting computable length bound and test triviality
in $E$. The resulting presentation uses the given generators, which need not
form a minimal generating set.
\end{proof}

\begin{lemma}\label{vfc:lem:certificate}
Suppose $H=\langle h_1,\ldots,h_k\mid R\rangle$ has a solvable word
problem. Let $E$ have a solvable word problem, and suppose a finite
presentation $\langle y_1,\ldots,y_k\mid S\rangle$ has been computed
for the subgroup $Q=\langle y_1,\ldots,y_k\rangle\leq E$.
The assignments $h_i\mapsto y_i$ define an isomorphism $H\to Q$ if
and only if
\begin{equation}\label{vfc:eq:checks}
R(y_1,\ldots,y_k)=1\text{ in }E,
\qquad
S(h_1,\ldots,h_k)=1\text{ in }H.
\end{equation}
Both conditions are decidable.
\end{lemma}

\begin{proof}
The first condition defines a homomorphism $H\to Q$ and the second
defines a homomorphism $Q\to H$. Their compositions fix the displayed
generators, so they are inverse. Necessity follows by applying an
isomorphism and its inverse to the relators.
\end{proof}

\begin{proposition}\label{vfc:prop:search}
There is a uniform algorithm on finite presentations of hyperbolic
groups satisfying the conclusion of Lemma~\ref{vfc:lem:structural} which
computes the data in Theorem~\ref{vfc:thm:main}(i)--(ii). In the second
case of that lemma the ambient group can be chosen hyperbolic and the
image quasiconvex. In the virtually free case the construction may give
$H=F_d\hookrightarrow F_d\times\Z$.
\end{proposition}

\begin{proof}
First compute a word-problem algorithm for $G$. Enumerate, for every
$n$, all assignments of the generators of $G$ to permutations in
$\operatorname{Sym}(n)$ satisfying its relators. Each assignment gives
a homomorphism $\rho$ to a finite group. Reidemeister--Schreier
rewriting computes a finite presentation
\[
H_\rho=\ker\rho=\langle h_1,\ldots,h_k\mid R_\rho\rangle
\]
and words for the $h_i$ in $G$. The finite image of $\rho$ gives the
index and a coset transversal. This enumerates every finite-index
normal subgroup. Equality of words in the $h_i$ is decided in $G$.

Run the following two searches in parallel.

In the first, enumerate finite sequences of Tietze transformations
starting from each presentation of $H_\rho$, looking for a presentation
of a finite-rank free group with no relators. Tietze transformations
are enumerated together with finite derivations justifying any added
consequence. A successful sequence gives explicit inverse maps
between $H_\rho$ and $F_d$. Return the embedding in $F_d\times\Z$
with zero $\Z$-coordinate. The rank-zero case gives $H_\rho=1$.

In the second, enumerate all ranks $d$ and all finite products of
Nielsen transformations of $F_d$. This enumerates $\Aut(F_d)$,
together with inverse automorphisms. Form $E_\alpha$ as in
\eqref{vfc:eq:ambient} and run hyperbolicity recognition on each such
presentation. Whenever hyperbolicity has been established for $E_\alpha$, enumerate,
for every $H_\rho$, all tuples $y_1,\ldots,y_k\in E_\alpha$.

Check $R_\rho(y_1,\ldots,y_k)=1$. If this holds, apply
Lemma~\ref{vfc:lem:qc}. When it produces a presentation with relators
$S$, check $S(h_1,\ldots,h_k)=1$ in $G$. A successful check gives
the required embedding by Lemma~\ref{vfc:lem:certificate}.

Interleave both searches and their subsidiary computations so that
each computation is eventually carried out for arbitrarily many steps.

If $G$ is virtually free, the normal core of a finite-index free
subgroup is free of finite rank. The first search eventually finds it.
Otherwise Lemma~\ref{vfc:lem:structural} supplies a finite-index normal
$H_\rho$ and a quasiconvex embedding into a hyperbolic $E_\alpha$.
Its generator images occur among the tuples in the second search.
The quasiconvexity procedure and both relator tests terminate for this
tuple, giving the required embedding.
\end{proof}

The first test in \eqref{vfc:eq:checks} defines a homomorphism to
$E_\alpha$. The presentation of its image and the second test give
the inverse homomorphism, and hence prove injectivity.

\subsection{A finite HNN presentation}\label{vfc:sec:hnn}

Identify the subgroup constructed in Proposition~\ref{vfc:prop:search}
with its image $H\leq E_\alpha$. Let $n_i=\pi(h_i)$. If every $n_i$ is
zero, then $H\leq F_d$, and Stallings folding computes a finite free
basis of $H$. We take $C=0$ and $K=H$.

Otherwise let
\[
q=\gcd(n_1,\ldots,n_k)>0,
\qquad
\chi=\frac{1}{q}\pi|_H:H\twoheadrightarrow\Z.
\]
B\'ezout coefficients give a word $s$ in the $h_i$ with $\pi(s)=q$.
Set $e_i=\chi(h_i)$ and $z_i=h_i s^{-e_i}$. Then
\[
H=\langle s,z_1,\ldots,z_k\rangle,
\qquad
K=\ker\chi=H\cap F_d.
\]
The subgroup theorem for free groups proves that $K$ is free.

Write $s=ct^q$ in the ambient normal form. Conjugation by $s$ induces
the explicitly computable automorphism
\begin{equation}\label{vfc:eq:beta}
\beta:F_d\longrightarrow F_d,
\qquad \beta(w)=c\alpha^q(w)c^{-1}.
\end{equation}
Its restriction preserves $K$ and satisfies $H=K\rtimes_\beta\langle s\rangle$.

\begin{proposition}\label{vfc:prop:hnn}
From the preceding data one can compute free groups $A,D\leq B\leq K$
of finite rank, an isomorphism $\theta:A\to D$, and a finite
presentation
\begin{equation}\label{vfc:eq:hnn}
H\cong
\left\langle b_1,\ldots,b_\ell,s\ \middle|\
sU_j(b)s^{-1}=V_j(b),\ 1\leq j\leq p\right\rangle,
\end{equation}
where the $b_i$ form a free basis of $B$, the $U_j$ form a free basis
of $A$, and $V_j=\theta(U_j)$ form a free basis of $D$. The
isomorphism with the original presentation of $H$ is computable in
both directions, and the HNN height is $\chi$.
\end{proposition}

\begin{proof}
Adjoin the generator $s$ with its known defining word and replace
$h_i$ by $z_i s^{e_i}$. This gives a finite presentation
\[
\mathcal P=\langle s,z_1,\ldots,z_k\mid R'\rangle
\]
in which each relator has exponent sum zero in $s$. Choose $N\geq1$
at least the absolute value of every prefix exponent sum in $s$ of
every relator in $R'$. For $-N\leq j\leq N$, compute the elements
\[
z_{i,j}=\beta^j(z_i)=s^jz_i s^{-j}\in F_d.
\]
Define
\begin{align*}
B&=\langle z_{i,j}:1\leq i\leq k,\ -N\leq j\leq N\rangle,\\
A&=\langle z_{i,j}:1\leq i\leq k,\ -N\leq j<N\rangle,\\
D&=\langle z_{i,j}:1\leq i\leq k,\ -N<j\leq N\rangle.
\end{align*}
Then $\beta(A)=D$. All these are subgroups of the explicitly given
finite-rank free group $F_d$. Stallings folding \cite{Stallings}
computes their free bases and rewrites the basis elements of $A$ and
their $\beta$-images as words $U_j,V_j$ in a free basis of $B$.
Retain expressions for the $b_i$ in the original generators
$z_{i,j}$. These expressions can also be found by enumeration and
free reduction.

Let $P$ denote the group on the right of \eqref{vfc:eq:hnn}. It is the
HNN extension of $B$ with associated subgroups $A,D$ and isomorphism
$\beta|_A$. The inclusion $B\leq H$ and the chosen element $s\in H$
define a homomorphism $P\to H$. It is surjective
because its image contains $s$ and all $z_{i,0}$.

In $P$, the defining relations imply
\[
s z_{i,j}s^{-1}=z_{i,j+1}\qquad(-N\leq j<N).
\]
It follows, for positive and negative $j$, that
$z_{i,j}=s^jz_{i,0}s^{-j}$ in $P$. Rewrite a relator of $R'$ by
collecting the occurrences of $s$. By the choice of $N$, the result
is a word in the $z_{i,j}$ just defined. This word is trivial in the
subgroup $B\leq F_d$, since the original relator is trivial
in the embedded copy of $H$. It is therefore trivial in $P$.
Consequently the assignments $s\mapsto s$ and $z_i\mapsto z_{i,0}$
define a homomorphism $H\to P$.

The two maps fix $s$ and all $z_{i,0}$. They also fix every
$z_{i,j}$ by the conjugation identities, and hence every generator
$b_i$ of $B$. Thus the maps are inverse. All the substitutions and
rewritings used above are effective. Finally, sending $s$ to $1$ and
$B$ to $0$ agrees with $\chi$ on the original generators.
\end{proof}

\begin{remark}
The base $B$ in \eqref{vfc:eq:hnn} has finite rank even when $K$ has
infinite rank. The isomorphism $\theta$ is between the two indicated
subgroups of $B$; it need not be an automorphism of $B$. The height kernel is free because
$H$ embeds in $E_\alpha$. This need not hold for an arbitrary HNN
extension of finite-rank free groups.
\end{remark}

\subsection{The free kernel and its automorphism}

\begin{proposition}\label{vfc:prop:kernel}
If $H\leq E_\alpha$ is quasiconvex and $\pi(H)\neq0$, the preceding
construction gives an effective free basis $\mathcal B$ of
$K=H\cap F_d$. Given $k\in K$ as a word in the ambient generators,
one can express $k$ as a word in $\mathcal B$. On every basis element
one can compute $\beta$ and $\beta^{-1}$ as words in $\mathcal B$.
\end{proposition}

\begin{proof}
Lemma~\ref{vfc:lem:qc} gives a membership algorithm for $H$ in
$E_\alpha$. For $w\in F_d$, membership in $K$ is exactly membership
in $H$, so it too is decidable.

Consider the right-coset Schreier graph of $K$ in $F_d$ with respect
to $a_1,\ldots,a_d$. Its vertices are $Kw$, and there is an edge
labelled $a_i$ from $Kw$ to $Kwa_i$. Equality of vertices is decidable,
since
\[
Kw=Kv\quad\Longleftrightarrow\quad wv^{-1}\in K.
\]
Represent each vertex by the shortest, and then lexicographically
least, reduced word in its coset. The representative of a given coset
is computable by a finite search: its length is at most that of the
given word. These representatives form a prefix-closed transversal,
and determine a computable spanning tree rooted at $K$.

For every unoriented edge outside this tree, choose its orientation
labelled by a positive basis letter. If the edge goes from $Ku$ to
$Kv$ with label $a_i$, and $u,v$ are the chosen representatives, put
\[
c_e=ua_iv^{-1}\in K.
\]
The Schreier basis theorem identifies
$\mathcal B=\{c_e:e\text{ is outside the spanning tree}\}$ as a
free basis of $K$. The chosen coset representatives and labelled edges can be enumerated,
and membership in the spanning tree is decidable. Thus the non-tree
edges give a decidable index set and an effective enumeration of
$\mathcal B$, whether or not this basis is finite.

A word representing an element of $K$ gives a finite closed path in
the Schreier graph. Reading its non-tree edges expresses the element
as a word in $\mathcal B$. Formula~\eqref{vfc:eq:beta} and its inverse
compute $\beta(c_e)$ and $\beta^{-1}(c_e)$ in $F_d$. They lie in $K$
because $s$ normalizes $K$, and the same rewriting expresses both in
$\mathcal B$.
\end{proof}

We obtain the effective presentation
\begin{equation}\label{vfc:eq:infinite}
H=\left\langle s,c_e\ (e\in J)\ \middle|\
sc_es^{-1}=W_e(\mathcal B)\ (e\in J)\right\rangle,
\end{equation}
where $J$ is the decidable index set above and $W_e$ is computed from
$e$. The isomorphisms from
\eqref{vfc:eq:hnn} and \eqref{vfc:eq:infinite} to the original
finite-index subgroup are explicit.

\begin{proof}[Proof of Theorem~\ref{vfc:thm:main}]
The assumed virtual compact specialness result and
\cite[Corollary~1.3]{KL} give the hypotheses of
Lemma~\ref{vfc:lem:structural}. Proposition~\ref{vfc:prop:search}
therefore constructs (i)--(ii). If the first
search succeeds, use the embedding $H=F_d\leq F_d\times\Z$ and
take $C=0$. For the second search, use the projection $\pi$ and the
greatest-common-divisor construction of Subsection~\ref{vfc:sec:hnn}. When its image
is zero, fold the generators of $H$ inside $F_d$. When its image is
nonzero, apply Propositions~\ref{vfc:prop:hnn} and~\ref{vfc:prop:kernel}.
These constructions give all the remaining data.
\end{proof}

\subsection{Finite-rank kernels}

If a finite-index subgroup with finite-rank free kernel is known to
exist, it can be found directly by the following procedure.

\begin{proposition}\label{vfc:prop:finite}
There is a uniform partial algorithm on finite presentations which
terminates whenever the input group has a finite-index subgroup
isomorphic to $F_e\rtimes_\gamma\Z$ with $e<\infty$. Upon
termination it returns a finite-index normal subgroup of this form,
its mapping-torus presentation, and explicit isomorphisms with its
Reidemeister--Schreier presentation.
\end{proposition}

\begin{proof}
Enumerate the kernels of all finite permutation representations and
their Reidemeister--Schreier presentations. Enumerate all ranks $e$
and products of Nielsen transformations $\gamma\in\Aut(F_e)$,
and the corresponding presentations \eqref{vfc:eq:ambient}. Search for
finite Tietze sequences joining the subgroup presentations to these
mapping-torus presentations, interleaving all searches.

It suffices to enumerate normal subgroups. Suppose $L\leq G$ has finite index and
$L=F_e\rtimes\Z$. Let $H$ be its normal core in $G$. Then
$H\cap F_e$ has finite index in $F_e$, and therefore has finite
rank. The image of $H$ in $\Z$ has finite nonzero index. The
resulting extension splits, so $H$ is itself a finite-rank
free-by-$\Z$ group. Its presentation eventually occurs in the
search and is joined by a finite Tietze sequence to a mapping-torus
presentation. A successful sequence provides the required maps.
\end{proof}

\begin{corollary}
Assume the result stated as Corollary~4.10 in \cite{OAISpecial}.
Let $F$ be a finite-rank
free group and let $\varphi:F\to F$ be an injective endomorphism.
If its ascending HNN extension is hyperbolic, one can compute a
finite-index normal subgroup with a finite-rank free-by-$\Z$ presentation.
\end{corollary}

\begin{proof}
By the cited result, the ascending HNN extension has a finite-index
subgroup with finite-rank free kernel. Apply
Proposition~\ref{vfc:prop:finite} to its finite presentation.
\end{proof}

The character obtained in this corollary need not be the restriction
of the original ascending HNN height.

\begin{example}
The finite-rank conclusion cannot hold for every hyperbolic
one-relator group. Let $G$ be the fundamental group of a closed
orientable surface of genus $g\geq2$. For every finite-index subgroup
$H\leq G$,
\[
\chi(H)=[G:H](2-2g)<0.
\]
A finite-rank free-by-$\Z$ group has a finite classifying space
given by the mapping torus of a homotopy equivalence of a finite
graph, and its Euler characteristic is zero. Thus no finite-index
subgroup of this surface group is finite-rank free-by-$\Z$.
Thus the free kernel in Theorem~\ref{vfc:thm:main} cannot in general
be required to have finite rank.
\end{example}

\begin{remark}
Termination follows from the existence theorems, without an a priori
bound on the finite index, the ambient rank, or the lengths of the
generator images. The argument gives no bound on the running time.
Proposition~\ref{vfc:prop:finite} is a semialgorithm; it does not decide
whether an arbitrary finitely presented group is virtually
finite-rank free-by-$\Z$.
\end{remark}

\section{Toral relative hyperbolicity of virtually free-by-cyclic groups}\label{sec:toral-vfc}

We prove that a finitely generated torsion-free virtually free-by-cyclic
group is toral relatively hyperbolic if and only if it is CSA. As in
Section~\ref{sec:effective-vfc}, the free kernel may have infinite rank.

A group is \emph{toral relatively hyperbolic} if it is torsion-free and
hyperbolic relative to a finite family of finitely generated abelian
subgroups. We allow the empty peripheral family and, for an abelian
group, the family consisting of the group itself. The Klein bottle group is
\[
K=\langle a,s\mid sas^{-1}=a^{-1}\rangle.
\]

\begin{theorem}\label{toral:thm:main}
Let $G$ be a finitely generated torsion-free virtually free-by-cyclic group.
The following are equivalent:
\begin{enumerate}
\item $G$ is toral relatively hyperbolic;
\item $G$ is CSA;
\item $G$ contains neither $F_2\times\Z$ nor $K$.
\end{enumerate}
When $G$ is nonabelian, its noncyclic peripheral groups may all be chosen
isomorphic to $\Z^2$.
\end{theorem}

We first prove the assertion for finite-rank free kernels, then pass
to arbitrary-rank kernels and finite extensions.

\subsection{Finite-rank free kernels}

\begin{lemma}\label{toral:lem:finite}
Let $F_r$ be a free group of finite rank and let
$J=F_r\rtimes_\alpha\langle t\rangle$.
If $J$ contains neither $F_2\times\Z$ nor $K$, then $J$ is toral relatively
hyperbolic, with noncyclic peripheral subgroups isomorphic to $\Z^2$.
\end{lemma}

\begin{proof}
Suppose a polynomial subgroup for the outer class of $\alpha$ has rank at
least two. By \cite[Lemma~1.5]{LevittGrowth}, there are $q\geq1$ and $g\in F_r$
such that
\[
\beta=\Inn(g)\circ\alpha^q,
\qquad \rank\Fix(\beta)\geq2.
\]
Here $\Inn(g)(x)=gxg^{-1}$ and $txt^{-1}=\alpha(x)$.
Choose noncommuting $x,y\in\Fix(\beta)$ and put $s=gt^q$.
Then $s$ centralizes $\langle x,y\rangle\cong F_2$.
Projection $J\to\Z$ shows that
$\langle s\rangle\cap\langle x,y\rangle=1$.
Consequently
\[
\langle x,y,s\rangle\cong F_2\times\Z,
\]
a contradiction. Thus every maximal polynomial subgroup is cyclic.

If $r\geq2$, $\alpha$ cannot itself be polynomially growing, since $F_r$
would then be a polynomial subgroup of rank at least two. The relative
hyperbolicity theorem for free-group automorphisms,
\cite[Theorem~3.5]{DL}, makes $J$ hyperbolic relative to the suspensions of
the finitely many maximal polynomial subgroups. A suspension of an infinite
cyclic polynomial subgroup has the form
\[
\Z\rtimes_{\varepsilon}\Z,\qquad \varepsilon\in\{1,-1\}.
\]
It is $\Z^2$ or $K$. The second possibility is excluded. If there are no
nontrivial polynomial subgroups, $J$ is hyperbolic.

For $r=0$, $J\cong\Z$. For $r=1$, the only possibilities are $\Z^2$ and
$K$, and the latter is excluded. This completes the proof.
\end{proof}

\subsection{Free kernels of arbitrary rank}

\begin{lemma}\label{toral:lem:infinite}
Let $H$ be finitely generated and free-by-cyclic. If $H$ contains neither
$F_2\times\Z$ nor $K$, then $H$ is toral relatively hyperbolic, with
noncyclic peripheral subgroups isomorphic to $\Z^2$.
\end{lemma}

\begin{proof}
The assertion is immediate if $H$ is free. Otherwise write
$H=F\rtimes_\alpha\Z$, allowing $F$ to have infinite rank.
By \cite[Theorem~1.1 and Remark~6.2]{LintonMT}, $H$ is hyperbolic relative
to a finite collection $\mathcal P$ whose members are finite-rank
free-by-cyclic groups. Thus the free kernels of the peripheral groups are finitely
generated. This also follows from
\cite[Corollary~4.10]{LintonMT}.

No $P\in\mathcal P$ contains $F_2\times\Z$ or $K$. By
Lemma~\ref{toral:lem:finite}, each $P$ is toral relatively hyperbolic
with noncyclic peripheral groups $\Z^2$. Transitivity of relative
hyperbolicity then gives the required peripheral structure for $H$.
If $\mathcal P=\{H\}$, apply the finite-rank case directly.
\end{proof}

The transitivity assertion follows from finite relative presentations.
Substitute presentations of the groups $P\in\mathcal P$ relative to
their abelian peripheral subgroups into a presentation of $H$ relative
to $\mathcal P$. The resulting relative presentation has a linear
isoperimetric inequality; see the proof of \cite[Theorem~3.5]{DL}.

\subsection{Passing to finite extensions}

\begin{lemma}\label{toral:lem:virtz2}
A torsion-free virtually $\Z^2$ group containing no Klein bottle subgroup
is isomorphic to $\Z^2$.
\end{lemma}

\begin{proof}
Let $A\cong\Z^2$ be a normal subgroup of finite index in $Q$.
Conjugation gives a finite-image homomorphism
$\rho:Q\to\operatorname{GL}(2,\Z)$.
For every $g\neq1$ there is $m\geq1$ with $1\neq g^m\in A$.
The element $g^m$ is fixed by conjugation by $g$, so $\rho(g)$ has
eigenvalue $1$.

A nonidentity finite-order matrix in $\operatorname{GL}(2,\Z)$ with
eigenvalue $1$ has other eigenvalue $-1$. If $\rho(g)$ is such a matrix,
choose $1\neq a\in A$ in its $-1$ eigenspace. Then $gag^{-1}=a^{-1}$.
Moreover, $\langle g\rangle\cap\langle a\rangle=1$: an equality
$g^i=a^j$ would, after conjugating by $g$, give $a^j=a^{-j}$, and
torsion-freeness gives $j=i=0$. Thus $\langle a,g\rangle\cong K$.
This is excluded, so $\rho$ is trivial.

Hence $A\leq Z(Q)$ and $Q/Z(Q)$ is finite. Schur's theorem makes $Q'$
finite; torsion-freeness gives $Q'=1$. Therefore $Q$ is a finitely generated
torsion-free abelian group, of rank two because $A$ has finite index.
\end{proof}

\begin{lemma}\label{toral:lem:finiteindex}
Suppose $G$ is finitely generated and torsion-free, and $H\leq G$ has
finite index. If $H$ is toral relatively hyperbolic with all noncyclic
peripheral subgroups isomorphic to $\Z^2$, and $G$ contains no $K$, then $G$ is toral relatively
hyperbolic with all noncyclic peripheral subgroups isomorphic to $\Z^2$.
\end{lemma}

\begin{proof}
Replace $H$ by its normal core. A finite-index subgroup of a toral relatively
hyperbolic group has the induced toral peripheral structure. If $H$ is
hyperbolic, so is $G$. If $H$ is virtually abelian, apply
Lemma~\ref{toral:lem:virtz2} (or the virtually cyclic case).

In the remaining case let $P_1,\ldots,P_k$ represent the conjugacy classes
of noncyclic maximal abelian subgroups of $H$. Conjugation by $G$ permutes
their $H$-conjugacy classes. Set
\[
Q_i=\operatorname{N}_G(P_i).
\]
Since $H$ is toral relatively hyperbolic, $\operatorname{N}_H(P_i)=P_i$. Thus
$Q_i\cap H=P_i$ and $[Q_i:P_i]<\infty$. Relative hyperbolicity passes
to $G$ with peripherals the $G$-conjugacy representatives among the $Q_i$.

This follows from the asymptotically tree-graded characterization of
relative hyperbolicity in \cite{Drutu}.
The inclusion $H\hookrightarrow G$ is a quasi-isometry. Since there are
finitely many cosets of $H$ and conjugation preserves the peripheral
family, the left cosets of the $Q_i$ correspond, at uniformly bounded
Hausdorff distance, to the left cosets of the $P_i$ in $H$.
The asymptotically tree-graded structure therefore passes to $G$ with
these peripheral cosets.

Finally each $Q_i$ is torsion-free and virtually $\Z^2$. By
Lemma~\ref{toral:lem:virtz2}, $Q_i\cong\Z^2$.
\end{proof}

\begin{proof}[Proof of Theorem~\ref{toral:thm:main}]
Assume (3). Choose a finite-index free-by-cyclic subgroup $H\leq G$.
It is finitely generated and inherits (3).
Lemma~\ref{toral:lem:infinite} applies to $H$, and
Lemma~\ref{toral:lem:finiteindex} applies to $H\leq G$. Thus (1) holds.

For a torsion-free group hyperbolic relative to abelian groups, every
nontrivial centralizer is either cyclic or contained in its unique
abelian peripheral conjugate. The almost malnormality of peripherals
and of maximal loxodromic elementary subgroups gives the CSA property.
Hence (1) implies (2).

If $G$ is CSA, it cannot contain $F_2\times\Z$, whose nontrivial central
element has nonabelian centralizer. It cannot contain $K$ either: if
$sas^{-1}=a^{-1}$, then $s$ normalizes the maximal abelian subgroup
containing $a$. Malnormality puts $s$ in that abelian subgroup, so
$a=a^{-1}$, contrary to torsion-freeness. Thus (2) implies (3).
\end{proof}

\begin{corollary}\label{toral:cor:onerelator}
Every finitely generated torsion-free CSA one-relator group which is
virtually free-by-cyclic is toral relatively hyperbolic. In particular,
this applies to
\[
G=\langle X\mid [u,v]\rangle
\]
whenever virtual free-by-cyclicity is established independently.
\end{corollary}

\begin{corollary}
Every finitely generated subgroup of a torsion-free CSA virtually
free-by-cyclic group is itself toral relatively hyperbolic.
\end{corollary}

\begin{proof}
If $L\leq G$ is finitely generated and $H\leq G$ is finite-index
free-by-cyclic, then $L\cap H$ has finite index in $L$. Its intersection
with the free kernel of $H$ is free and its image in the cyclic quotient
is cyclic. Hence $L$ is virtually free-by-cyclic. It also inherits
torsion-freeness and CSA, so Theorem~\ref{toral:thm:main} applies.
\end{proof}

\begin{remark}
The subgroups in the last corollary need not be relatively quasiconvex
in $G$. For example, a free fibre of a finite-rank free-by-cyclic group
can be distorted.
\end{remark}

\subsection{CSA one-relator groups}

The results above do not establish toral relative hyperbolicity for
every torsion-free CSA one-relator group containing $\Z^2$, even when
the relator is a commutator. Section~\ref{sec:effective-vfc} assumes
hyperbolicity and the virtual compact specialness theorem, so it does
not give virtual free-by-cyclicity in this case.

We recall a reduction using one-relator towers. By
\cite[Propositions~4.3 and~4.6, Lemmas~4.4 and~5.2]{LintonOR}, a factored
one-relator tower is stable down to a primitive extension core, unless it
is stable throughout or encounters a powered core. A powered core is
noncyclic and generated by $a,b$ satisfying $a^m=b^n\neq1$.
It cannot occur in a torsion-free CSA group: $a,b$ centralize the common
power, so they commute, and the resulting two-generated torsion-free
abelian group has rank at most one.

It remains to prove toral relative hyperbolicity for the two-generator
primitive extension cores not satisfying the hypotheses above, and to
prove relative quasiconvexity of the Magnus subgroups used in the
stable part of the tower. The condition on segment stabilizers alone
does not imply this relative quasiconvexity.

\section{Algorithmic cases for CSA one-relator groups}\label{sec:csa-algorithms}

Let $G=\langle X\mid r\rangle$ be a torsion-free one-relator group.
We consider the property
\[
 \PP(G):\qquad G\text{ is CSA and contains a subgroup isomorphic to }\Z^2.
\]
An algorithm on a specified class is required to terminate on every
presentation defining a group in that class. It need not decide
membership in the class. A semialgorithm for a property terminates
precisely when that property holds.

\subsection{A recognition algorithm}

\begin{proposition}\label{cases:prop:algorithm}
Let $\mathcal C$ be a class of torsion-free one-relator groups such that every
CSA group in $\mathcal C$ is toral relatively hyperbolic. Then there is a uniform
algorithm which, on presentations of groups in $\mathcal C$,
decides the CSA property and, when the answer is affirmative, decides whether
$\Z^2\leq G$. In particular, it decides $\PP(G)$.
\end{proposition}

\begin{proof}
Use the uniform solution of the word problem for one-relator groups.
First check whether the generators commute pairwise. If they do, $G$ is
abelian, hence CSA, and Smith normal form computes its rank and decides
whether it contains $\Z^2$.

Otherwise run the following two searches in parallel.
\begin{enumerate}[label=\textup{(\roman*)},leftmargin=*,itemsep=0.35em]
\item Enumerate word triples $(a,b,g)$ and look for
\begin{equation}\label{cases:eq:failure}
 a\ne1,\quad b\ne1,\quad [a,b]=[a,gbg^{-1}]=1,
 \quad [a,g]\ne1.
\end{equation}
Such a triple proves that $G$ is not CSA, so $\PP(G)$ is false.
\item Run Dahmani's recognition algorithm for relative hyperbolicity with
abelian peripheral groups \cite{Dahmani}. When it returns an exact relative
presentation, compute the ranks of the peripheral groups. Then $G$ is CSA, and $\PP(G)$ holds if and only if some peripheral
subgroup has rank at least two.
\end{enumerate}

A group is not CSA if and only if a triple satisfying
\eqref{cases:eq:failure} exists.
Indeed, in a CSA group, $b$ and $gbg^{-1}$ belong to the maximal abelian
subgroup containing $a$; malnormality forces $g$ into that subgroup.
Conversely, if no such triple exists, the centralizer of every $b\ne1$ is
abelian: apply the condition to $a,g\in C_G(b)$. These abelian centralizers
are malnormal by applying the same condition to a nontrivial element in an
intersection with a conjugate. Thus $G$ is CSA.

Search (i) terminates if $G$ is not CSA. If $G$ is CSA and belongs to
$\mathcal C$, search (ii) terminates by the hypothesis. A torsion-free group hyperbolic relative to abelian groups is CSA
\cite[Lemma~10.1]{DGroves}, and every $\Z^2$ subgroup is parabolic.
The exact peripheral presentations determine the ranks of these
subgroups. Thus both answers are correct.
\end{proof}

The procedure decides both properties on $\mathcal C$. As in
Section~\ref{sec:recognition}, the proof gives no bound on its running time.

\subsection{Applications}

We apply Proposition~\ref{cases:prop:algorithm} to the following classes.
Except in the first case, CSA is a property to be decided.
The conclusions in Cases~\ref{cases:cyclic}--\ref{cases:exceptional}
and the finite-hierarchy application assume the combination statement
of \cite{Stable}, with the hypotheses given in Case~\ref{cases:stable}.

\begin{enumerate}[label=\textbf{Case \arabic*.},ref=\arabic*,leftmargin=*,itemsep=0.75em]
\item\label{cases:toral} \textbf{Toral relatively hyperbolic groups.}
Such a group is CSA. Dahmani's algorithm computes a peripheral
structure, and the ranks of the peripheral groups decide $\PP(G)$.

\item\label{cases:vfc} \textbf{Virtually free-by-cyclic groups.}
Suppose $G$ has a finite-index subgroup which is free or has the form
\[
 F\rtimes_\alpha\Z,
\]
where $F$ is a free group, possibly of infinite rank. Theorem~\ref{toral:thm:main}
gives, for finitely generated torsion-free groups in this class,
\[
 G\text{ is CSA}\quad\Longleftrightarrow\quad
 G\text{ is toral relatively hyperbolic}.
\]
Proposition~\ref{cases:prop:algorithm} therefore decides both CSA and
$\PP(G)$ in this class.

\item\label{cases:ascending} \textbf{Ascending HNN extensions of finite-rank free groups.}
Suppose
\[
 G\cong\langle F,t\mid txt^{-1}=\phi(x)\ (x\in F)\rangle,
 \qquad \phi:F\hookrightarrow F,
\]
with $F$ finitely generated free. Surjectivity is not required. If $G$ is
CSA, it contains no $\BS(1,d)$ for $d\geq2$. Mutanguha's relative expansion
theorem gives free-by-cyclic peripheral groups. By Case~\ref{cases:vfc},
these are toral relatively hyperbolic, so transitivity gives abelian
peripheral groups \cite[Theorem~4.7]{Mutanguha}. The automorphism case is
already included in Case~\ref{cases:vfc}.

\item\label{cases:partial} \textbf{Partial ascending extensions with a free-factor domain.}
Suppose
\[
 G\cong\langle F,t\mid tat^{-1}=\psi(a)\ (a\in A)\rangle,
 \qquad F=A*B,
\]
where $F$ is finitely generated free and $\psi:A\hookrightarrow F$ is
injective. Chong--Wise embed $G$ as an undistorted peripheral subgroup of
a finite-rank free ascending HNN extension $L$ \cite{CW}. If $G$ is CSA,
so is $L$; Case~\ref{cases:ascending} makes $L$ toral relatively hyperbolic. Undistortion and
Hruska's theorem then make $G$ relatively quasiconvex for that toral
structure, hence toral relatively hyperbolic \cite{Hruska}.
The cited embedding theorem requires $A$ to be a free factor of $F$.

\item\label{cases:subgroups} \textbf{Finitely generated subgroups of free-group mapping tori.}
More generally, suppose $G$ embeds in the mapping torus of an injective
endomorphism of a free group of arbitrary rank. Linton's Theorem~1.1 makes
each finitely generated nonfree such subgroup relatively hyperbolic with
peripherals which are mapping tori of finite-rank free groups
\cite{LintonMT}. If $G$ is CSA, these peripherals inherit CSA; apply Case~\ref{cases:ascending}
and refine the peripheral structure. Free subgroups are handled directly.
The algorithm requires only the existence of this embedding; it does
not require the embedding as part of the input.

\item\label{cases:cyclic} \textbf{Finite graphs of free groups with cyclic edge groups.}
Suppose $G$ admits a finite graph-of-groups decomposition with finitely
generated free vertex groups and cyclic or trivial edge groups.
This includes cyclically pinched and conjugacy-pinched one-relator
presentations. It is a special case of Case~\ref{cases:stable}: the stability constant is
$N=1$, and finitely generated subgroups of free groups are quasiconvex.

\item\label{cases:stable} \textbf{Stable graphs of toral relatively hyperbolic groups.}
Suppose $G$ admits a finite graph-of-groups decomposition such that:
\begin{itemize}[leftmargin=1.4em,itemsep=0.15em]
\item vertex groups are finitely generated and toral relatively hyperbolic;
\item edge groups are finitely generated free and relatively quasiconvex
in both endpoint vertex groups;
\item for some $N\geq1$, every $N$-edge segment in the Bass--Serre tree has
cyclic or trivial pointwise stabilizer.
\end{itemize}
The stable combination statement in \cite{Stable} gives
CSA $\Rightarrow$ toral relative hyperbolicity under these hypotheses.
The algorithm applies when the displayed hypotheses hold. We do not
require an algorithm deciding whether an arbitrary one-relator group
admits such a decomposition.

\item\label{cases:overlapping} \textbf{Overlapping commutator relators.}
Let $A,B,C$ be disjoint finite alphabets, with
\[
 x\in F(A\cup B)\setminus F(B),\qquad
 y\in F(B\cup C)\setminus F(B),\qquad
 1\ne z\in F(B).
\]
The groups
\[
 G=\langle A,B,C\mid[xy,z]=1\rangle
\]
satisfy the hypotheses of Case~\ref{cases:stable}. Rewrite the relation
as $x^{-1}zx=yzy^{-1}$.
The resulting amalgam has free vertex groups and a free edge group
isomorphic to $F(B)*\Z$. The two-edge segment stabilizers are cyclic or
trivial by the Magnus conjugate-intersection theorem. Thus Case~\ref{cases:stable} applies
with $N=2$. The displayed word conditions are decidable in the ambient
free group. This argument uses the displayed conditions on $x,y,z$; it does not
apply to an arbitrary commutator relator.

\item\label{cases:exceptional} \textbf{Exceptional Magnus intersections.}
Suppose the presentation has Magnus subgroups generated by $A\cup B$ and
$B\cup C$, with $A,B,C$ disjoint, such that
\[
 \langle A,B\rangle\cap\langle B,C\rangle
 \supsetneq \langle B\rangle.
\]
Linton's dichotomy gives either a $2$-free group or a two-generator
one-relator generalized Baumslag--Solitar subgroup $H$
\cite[Theorem~3.17 and Remark~3.18]{LintonOR}. The $2$-free case is
hyperbolic. In the other case, CSA forces the nonfree group $H$ to be
$\Z^2$: a CSA generalized Baumslag--Solitar group is abelian, since all
its cyclic vertex groups are commensurated. The power-equality alternative
in Howie's description is impossible for a generating pair of $\Z^2$;
the remaining alternative reduces to Case~\ref{cases:overlapping}. Thus CSA again implies
toral relative hyperbolicity. This application requires the displayed intersection to be exceptional.
\end{enumerate}

\paragraph{Finite hierarchies.}
The same conclusion holds for a finite hierarchy whose terminal groups
satisfy one of Cases~\ref{cases:toral}--\ref{cases:subgroups} and whose successive splittings satisfy the
free-edge, relative-quasiconvexity and segment-stabilizer assumptions of
Case~\ref{cases:stable}. If the total group is CSA, all vertex groups inherit CSA, so
induction gives toral relative hyperbolicity. The existence of such a hierarchy is a hypothesis of this application.

\subsection{Explicit examples and recognition of presentations}

\paragraph{A family of CSA groups containing $\Z^2$.}
For every integer $m\geq2$, the family
\[
 G_m=\langle b,t\mid[t,b^{-1}tbt^{-1}b^{-m}]=1\rangle
\]
satisfies $\PP(G_m)$. Introducing $a=b^{-1}tbt^{-1}b^{-m}$ gives
\[
 G_m\cong\langle a,b,t\mid[t,a]=1,\quad tbt^{-1}=bab^m\rangle.
\]
The free-group calculation and the relative expansion theorem give
relative hyperbolicity with the single peripheral subgroup
$\langle a,t\rangle\cong\Z^2$. Hence $G_m$ is CSA. The groups in this
family have noncyclic stabilizers of arbitrarily long segments in the
tree of the ascending HNN extension. Thus the ascending HNN argument
applies, whereas the segment-stabilizer hypothesis of
Case~\ref{cases:stable} does not hold.

For the standard Baumslag--Solitar presentations
$\BS(m,n)=\langle a,t\mid ta^mt^{-1}=a^n\rangle$, with $m,n\ne0$, the
answer is positive exactly when $(m,n)=(1,1)$ or $(-1,-1)$; both give
$\Z^2$. The other cases are not CSA \cite{GKM}.

\paragraph{Recognition from presentations.}
For finite-rank free-by-cyclic, ascending, and partial ascending
presentations with free-factor domain, free-group algorithms verify
injectivity and the required basis conditions. Finite graphs of free
groups with injective cyclic edge maps can also be enumerated. A finite
sequence of Tietze transformations verifies that one of these
presentations defines the original group. One can similarly enumerate
finite-index subgroups together with finite-rank free-by-cyclic
presentations. Applying Proposition~\ref{cases:prop:algorithm} to these
presentations gives semialgorithms for $\PP(G)$ within these classes.

This enumeration need not decide whether a group belongs to one of the
classes. For infinite-rank kernels in Case~\ref{cases:vfc}, unspecified
embeddings in Case~\ref{cases:subgroups}, and the geometric hypotheses of
Case~\ref{cases:stable}, the decision procedures apply under the stated
assumptions.

\subsection{Arbitrary torsion-free one-relator groups}

The following procedures apply to every torsion-free one-relator group:
\begin{itemize}[leftmargin=*,itemsep=0.35em]
\item The search for a triple satisfying \eqref{cases:eq:failure}
terminates precisely when $G$ is not CSA.
\item The hyperbolicity-detection procedure terminates on hyperbolic
groups. Such a group contains no $\Z^2$, so $\PP(G)$ is false.
\item Finitely many commutator tests decide whether $G$ is abelian;
its rank then decides $\PP(G)$.
\end{itemize}

Whenever the relative-hyperbolicity procedure produces an abelian
peripheral subgroup of rank at least two, $\PP(G)$ holds. To prove that
this procedure terminates for every group satisfying $\PP(G)$, one
would need to establish
\[
 \bigl(G\text{ torsion-free one-relator and CSA},\quad \Z^2\leq G\bigr)
 \ \Longrightarrow\ G\text{ toral relatively hyperbolic}
\]
in general. The results above do not prove this implication for every
commutator presentation $\langle X\mid[u,v]\rangle$.

The forbidden-subgroup criterion in \cite{GKM} gives a semialgorithm
for failure of CSA. A finite unsuccessful search does not prove CSA.
The decision algorithms above use, in addition, the implication from
CSA to toral relative hyperbolicity in each stated class.

\paragraph{Acknowledgment.}
The development of algorithms from the results stated in
\cite{OAI,OAISpecial} and the editing of this paper were carried out
with the assistance of ChatGPT~6~Pro. Our original model-theoretic results are
due to the authors and originate in our earlier manuscript, written
in 2020; they are independent of both cited results.


\begin{thebibliography}{99}
\raggedright

\bibitem{BMR}
G. Baumslag, A. Myasnikov, and V. Remeslennikov,
\emph{Algebraic geometry over groups I. Algebraic sets and ideal theory},
J. Algebra \textbf{219} (1999), 16--79.

\bibitem{Brown}
K.~S.~Brown,
\emph{Cohomology of Groups},
Graduate Texts in Mathematics, vol.~87, Springer-Verlag, New York, 1982.

\bibitem{BKM}
I. Bumagin, O. Kharlampovich, and A. Miasnikov,
\emph{Isomorphism problem for fully residually free groups},
J. Pure Appl. Algebra \textbf{208} (2007), no.~3, 961--977.

\bibitem{BK}
J.~O.~Button and R.~P.~Kropholler,
\emph{Nonhyperbolic free-by-cyclic and one-relator groups},
New York J. Math. \textbf{22} (2016), 755--774.
\href{https://nyjm.albany.edu/j/2016/22-35p.pdf}{Full text}.

\bibitem{CW}
H.~K.~Chong and D.~T.~Wise, \emph{Embedding partial HNN extensions in
ascending HNN extensions},
\href{https://arxiv.org/abs/2408.00453}{arXiv:2408.00453},
Theorem~1.1 and Corollaries~5.12--5.13.

\bibitem{CHR}
L. Ciobanu, D. Holt, and S. Rees,
\emph{Equations in groups that are virtually direct products},
J. Algebra \textbf{545} (2020), 88--99.

\bibitem{Dahmani}
F.~Dahmani, \emph{Finding relative hyperbolic structures},
Bull. London Math. Soc. \textbf{40} (2008), no.~3, 395--404.
\href{https://doi.org/10.1112/blms/bdn018}{doi:10.1112/blms/bdn018}.

\bibitem{DGroves}
F.~Dahmani and D.~Groves, \emph{The isomorphism problem for toral relatively
hyperbolic groups}, Publ. Math. IHES \textbf{107} (2008), 211--290,
\href{https://doi.org/10.1007/s10240-008-0014-3}{doi:10.1007/s10240-008-0014-3}.

\bibitem{DGEquations}
F. Dahmani and V. Guirardel,
\emph{Foliations for solving equations in groups: free, virtually free,
and hyperbolic groups}, J. Topol. \textbf{3} (2010), 343--404;
\href{https://doi.org/10.1112/jtopol/jtq010}{doi:10.1112/jtopol/jtq010}.

\bibitem{DG}
F. Dahmani and V. Guirardel,
\emph{The isomorphism problem for all hyperbolic groups},
Geom. Funct. Anal. \textbf{21} (2011), 223--300;
\href{https://doi.org/10.1007/s00039-011-0120-0}{doi:10.1007/s00039-011-0120-0}.

\bibitem{DL}
F.~Dahmani and R.~Li,
\emph{Relative hyperbolicity for automorphisms of free products and free groups},
J. Topol. Anal. \textbf{14} (2022), no.~1, 55--92.
\href{https://doi.org/10.1142/S1793525321500011}{doi:10.1142/S1793525321500011}.

\bibitem{DMBS}
E. Daniyarova and A. Myasnikov,
\emph{Groups elementarily equivalent to metabelian Baumslag--Solitar
groups and regular bi-interpretability}, Ann. Pure Appl. Logic
\textbf{177} (2026), no.~5, article 103695;
\href{https://doi.org/10.1016/j.apal.2025.103695}{doi:10.1016/j.apal.2025.103695}.

\bibitem{deLaNuez}
J. de la Nuez Gonz\'alez,
\emph{On expansions of non-abelian free groups by cosets of a finite
index subgroup}, doctoral dissertation, University of M\"unster, 2016.
\href{https://noah.nrw/ulbmshsnoah/content/titleinfo/4295670}{University repository}.

\bibitem{Drutu}
C.~Dru\c{t}u, \emph{Relatively hyperbolic groups: geometry and
quasi-isometric invariance},
Comment. Math. Helv. \textbf{84} (2009), no.~3, 503--546.
\href{https://doi.org/10.4171/CMH/171}{doi:10.4171/CMH/171}.

\bibitem{FV}
S. Feferman and R. L. Vaught,
\emph{The first order properties of products of algebraic systems},
Fund. Math. \textbf{47} (1959), 57--103;
\href{https://doi.org/10.4064/fm-47-1-57-103}{doi:10.4064/fm-47-1-57-103}.

\bibitem{GKL}
G. Gardam, D. Kielak, and A. D. Logan,
\emph{Algebraically hyperbolic groups}, Groups Geom. Dyn. (2025), published online first;
\href{https://doi.org/10.4171/GGD/907}{doi:10.4171/GGD/907}.

\bibitem{GardamWoodhouse}
G. Gardam and D. J. Woodhouse,
\emph{The geometry of one-relator groups satisfying a polynomial
isoperimetric inequality},
Proc. Amer. Math. Soc. \textbf{147} (2019), no.~1, 125--129.
\href{https://doi.org/10.1090/proc/14238}{doi:10.1090/proc/14238}.

\bibitem{Gildenhuys}
D.~Gildenhuys,
\emph{Classification of soluble groups of cohomological dimension two},
Math. Z. \textbf{166} (1979), no.~1, 21--25.
\url{https://doi.org/10.1007/BF01173844}.

\bibitem{GKM}
D. Gildenhuys, O. Kharlampovich, and A. Myasnikov,
\emph{CSA-groups and separated free constructions},
Bull. Austral. Math. Soc. \textbf{52} (1995), 63--84.

\bibitem{GS}
D.~Gildenhuys and R.~Strebel,
\emph{On the cohomological dimension of soluble groups},
Canad. Math. Bull. \textbf{24} (1981), no.~4, 385--392.
\url{https://doi.org/10.4153/CMB-1981-060-x}.

\bibitem{GrayLevine}
R. D. Gray and A. Levine,
\emph{The Diophantine problem in the Baumslag--Gersten group},
2026 preprint;
\href{https://arxiv.org/abs/2607.24189}{arXiv:2607.24189}, Theorem 1.1.

\bibitem{GWLimit}
D. Groves and H. Wilton,
\emph{Enumerating limit groups}, Groups Geom. Dyn. \textbf{3} (2009),
389--399;
\href{https://doi.org/10.4171/GGD/63}{doi:10.4171/GGD/63}.

\bibitem{GuyotStalder}
L. Guyot and Y. Stalder,
\emph{Limits of Baumslag--Solitar groups and dimension estimates
in the space of marked groups}, Groups Geom. Dyn. \textbf{6} (2012),
533--577;
\href{https://doi.org/10.4171/GGD/165}{doi:10.4171/GGD/165}, Proposition 5.18.

\bibitem{Hruska}
G.~C.~Hruska, \emph{Relative hyperbolicity and relative quasiconvexity for
countable groups}, Algebr. Geom. Topol. \textbf{10} (2010), 1807--1856,
Theorems~1.5 and~9.1.

\bibitem{Kap}
I.~Kapovich,
\emph{Detecting quasiconvexity: algorithmic aspects},
in \emph{Geometric and Computational Perspectives on Infinite Groups},
DIMACS Series in Discrete Mathematics and Theoretical Computer Science
\textbf{25}, American Mathematical Society, 1996, 91--99.

\bibitem{KMHyp}
O. Kharlampovich and A. Myasnikov,
\emph{Decidability of the elementary theory of a torsion-free hyperbolic group};
\href{https://arxiv.org/abs/1303.0760v8}{arXiv:1303.0760v8}, Theorems 1--3.

\bibitem{KMFree}
O. Kharlampovich and A. Myasnikov,
\emph{Elementary theory of free non-abelian groups},
J. Algebra \textbf{302} (2006), 451--552.

\bibitem{KMS} O. Kharlampovich, A.  Miasnikov and M. Sohrabi, Rich groups, weak second-order logic and applications, in Groups and Model Theory, Gagta Book 2, De Gruyter, 2020.

\bibitem{KMW}
O.~Kharlampovich, A.~Miasnikov and P.~Weil,
\emph{Stallings graphs for quasi-convex subgroups},
Journal of Algebra \textbf{488} (2017), 442--483.
\href{https://doi.org/10.1016/j.jalgebra.2017.05.037}{doi:10.1016/j.jalgebra.2017.05.037}.

\bibitem{Khelif}
A. Khelif,
\emph{Bi-interpr\'etabilit\'e et structures QFA: \'etude de groupes
r\'esolubles et des anneaux commutatifs},
C. R. Math. Acad. Sci. Paris \textbf{345} (2007), 59--61;
\href{https://doi.org/10.1016/j.crma.2007.06.003}{doi:10.1016/j.crma.2007.06.003}.

\bibitem{KL}
D.~Kielak and M.~Linton,
\emph{Virtually free-by-cyclic groups},
Geometric and Functional Analysis \textbf{34} (2024), 1580--1608.
\href{https://doi.org/10.1007/s00039-024-00687-6}{doi:10.1007/s00039-024-00687-6}.

\bibitem{Kropholler}
P.~H.~Kropholler,
\emph{Baumslag--Solitar groups and some other groups of cohomological
dimension two},
Comment. Math. Helv. \textbf{65} (1990), no.~4, 547--558.
\href{https://doi.org/10.1007/BF02566625}{doi:10.1007/BF02566625}.

\bibitem{KLL}
P.~H.~Kropholler, P.~A.~Linnell, and W.~L\"uck,
\emph{Groups of small homological dimension and the Atiyah conjecture},
in \emph{Geometric and Cohomological Methods in Group Theory},
London Math. Soc. Lecture Note Ser., vol.~358,
Cambridge University Press, 2009, pp.~272--277.

\bibitem{LevittGrowth}
G.~Levitt, \emph{Counting growth types of automorphisms of free groups},
Geom. Funct. Anal. \textbf{19} (2009), no.~4, 1119--1146.
\href{https://doi.org/10.1007/s00039-009-0016-4}{doi:10.1007/s00039-009-0016-4}.

\bibitem{Levitt}
G.~Levitt,
\emph{Quotients and subgroups of Baumslag--Solitar groups},
J. Group Theory \textbf{18} (2015), no.~1, 1--43.
\href{https://doi.org/10.1515/jgth-2014-0028}{doi:10.1515/jgth-2014-0028}.

\bibitem{LintonMT}
M.~Linton, \emph{The geometry of subgroups of mapping tori of free groups},
\href{https://arxiv.org/abs/2510.03145}{arXiv:2510.03145},
Theorem~1.1, Corollary~4.10, and Remark~6.2.

\bibitem{LintonOR}
M.~Linton, \emph{Hyperbolic one-relator groups},
Canad. J. Math. \textbf{78} (2026), no.~1, 35--61.
\href{https://doi.org/10.4153/S0008414X24000427}{doi:10.4153/S0008414X24000427}.

\bibitem{Linton}
M. Linton, \emph{One-relator hierarchies},
Duke Math. J. \textbf{174} (2025), no. 4, 747--802;
\href{https://doi.org/10.1215/00127094-2024-0040}{doi:10.1215/00127094-2024-0040}.

\bibitem{Lyndon}
R.~C.~Lyndon,
\emph{Cohomology theory of groups with a single defining relation},
Ann. of Math. (2) \textbf{52} (1950), 650--665.
\href{https://doi.org/10.2307/1969440}{doi:10.2307/1969440}.

\bibitem{LS}
R. C. Lyndon and P. E. Schupp,
\emph{Combinatorial Group Theory}, Springer, 1977.

\bibitem{MU}
R. Mandel and A. Ushakov,
\emph{Quadratic equations in metabelian Baumslag--Solitar groups},
Internat. J. Algebra Comput. \textbf{33} (2023), 1195--1216.

\bibitem{MZ}
A. Minasyan and P. Zalesskii,
\emph{One-relator groups with torsion are conjugacy separable},
J. Algebra \textbf{382} (2013), 39--45.

\bibitem{Mutanguha}
J.~P.~Mutanguha, \emph{Relative expansions for free group endomorphisms},
October~9, 2024,
\href{https://mutanguha.com/pdfs/survey.pdf}{author's survey}, Theorem~4.7.

\bibitem{Newman}
B. B. Newman, \emph{Some results on one-relator groups},
Bull. Amer. Math. Soc. \textbf{74} (1968), 568--571;
\href{https://doi.org/10.1090/S0002-9904-1968-12012-9}{doi:10.1090/S0002-9904-1968-12012-9}.

\bibitem{NybergBG}
C.-F. Nyberg-Brodda,
\emph{Undecidability of the Diophantine problem for one-relator groups
and one-relation monoids}, 2026 preprint;
\href{https://arxiv.org/abs/2608.01983}{arXiv:2608.01983}, Theorem A.

\bibitem{OAI}
OpenAI, \emph{Baumslag--Solitar-free one-relator groups are hyperbolic},
preprint, 25 September 2026.

\bibitem{OAISpecial}
OpenAI,
\emph{Virtual compact specialness of hyperbolic one-relator groups},
preprint, September 25, 2026. Theorem~1.1, Corollary~4.10 and
Corollary~18.1. The theorem numbers refer to the supplied version.

\bibitem{Papasoglu}
P. Papasoglu, \emph{An algorithm detecting hyperbolicity}, in
\emph{Geometric and Computational Perspectives on Infinite Groups},
DIMACS Ser. Discrete Math. Theoret. Comput. Sci., vol. 25,
Amer. Math. Soc., Providence, RI, 1996, pp. 193--200.

\bibitem{JR}
J.~Robinson,
\emph{Definability and decision problems in arithmetic},
Journal of Symbolic Logic \textbf{14} (1949), no.~2, 98--114.

\bibitem{SelaVII}
Z. Sela, \emph{Diophantine geometry over groups VII: the elementary
theory of a hyperbolic group}, Proc. Lond. Math. Soc. \textbf{99} (2009),
217--273;
\href{https://doi.org/10.1112/plms/pdn052}{doi:10.1112/plms/pdn052}.

\bibitem{SelaHopf}
Z. Sela, \emph{Endomorphisms of hyperbolic groups I: the Hopf property},
Topology \textbf{38} (1999), 301--321;
\href{https://doi.org/10.1016/S0040-9383(98)00015-9}{doi:10.1016/S0040-9383(98)00015-9}.

\bibitem{Serre}
J.-P.~Serre,
\emph{Trees},
translated from the French by John Stillwell,
Springer-Verlag, Berlin--New York, 1980.

\bibitem{Stable}
\emph{Stable toral combination and peripheral completion},
companion manuscript, October~14, 2026.
The combination statement is used under the hypotheses of
Case~\ref{cases:stable}.

\bibitem{Stallings}
J.~R.~Stallings,
\emph{Topology of finite graphs},
Inventiones Mathematicae \textbf{71} (1983), 551--565.
\href{https://doi.org/10.1007/BF02095993}{doi:10.1007/BF02095993}.

\bibitem{TMR}
A.~Tarski, A.~Mostowski, and R.~M.~Robinson,
\emph{Undecidable Theories},
Studies in Logic and the Foundations of Mathematics,
North-Holland Publishing Company, Amsterdam, 1953.

\bibitem{vdDG}
L.~van den Dries and A.~G\"unayd\i n,
\emph{The fields of real and complex numbers with a small multiplicative
group}, Proc. London Math. Soc. (3) \textbf{93} (2006), no.~1, 43--81.
\url{https://web.bogazici.edu.tr/ayhan.gunaydin/Real.pdf}.

\bibitem{WangZhang}
K. Wang and Q. Zhang,
\emph{Commutator relators of one-relator groups do not force Hopficity,
residual finiteness, or automaticity}, 2026 preprint;
\href{https://arxiv.org/abs/2607.21493}{arXiv:2607.21493}, Theorem 1.3.

\bibitem{WiltonSurface}
H. Wilton, \emph{Surface groups among cubulated hyperbolic and one-relator
groups}, Acta Math., to appear;
\href{https://arxiv.org/abs/2406.02121v3}{arXiv:2406.02121v3}, Theorem D.

\bibitem{Wise}
D. T. Wise, \emph{The Structure of Groups with a Quasiconvex Hierarchy},
Annals of Mathematics Studies 209, Princeton University Press, 2021.

\end{thebibliography}
\end{document}